\documentclass[review,onefignum,onetabnum]{siamart250211}

\usepackage{lipsum}
\usepackage{amsfonts}
\usepackage{graphicx}
\usepackage{epstopdf}
\usepackage{algorithmic}
\usepackage{amssymb}
\usepackage[english]{babel}
\usepackage{array}
\usepackage{adjustbox}
\usepackage{relsize}
\usepackage{spverbatim}
\usepackage{listings}
\usepackage{mathrsfs}
\usepackage{derivative}
\usepackage{mathtools}
\usepackage{float}
\usepackage[section]{placeins}

\usepackage{hyperref}
\usepackage{cleveref}

\usepackage{etoolbox}

\makeatletter
\newcommand{\crefcomma}[1]{%
  \begingroup
    \def\crefcomma@sep{}%
    \forcsvlist{\crefcomma@do}{#1}%
  \endgroup
}
\newcommand{\crefcomma@do}[1]{%
  \ifx\crefcomma@sep\@empty\else,~\fi
  \cref{#1}%
  \def\crefcomma@sep{,}%
}
\newcommand{\Crefcomma}[1]{%
  \begingroup
    \def\crefcomma@sep{}%
    \forcsvlist{\Crefcomma@do}{#1}%
  \endgroup
}
\newcommand{\Crefcomma@do}[1]{%
  \ifx\crefcomma@sep\@empty\else,~\fi
  \Cref{#1}%
  \def\crefcomma@sep{,}%
}
\makeatother

\DeclareMathAlphabet{\mathpzc}{OT1}{pzc}{m}{it}

\newcolumntype{L}{>{$}l<{$}}

\ifpdf
  \DeclareGraphicsExtensions{.eps,.pdf,.png,.jpg}
\else
  \DeclareGraphicsExtensions{.eps}
\fi

\newsiamremark{hypothesis}{Hypothesis}
\crefname{hypothesis}{Hypothesis}{Hypotheses}
\newsiamthm{claim}{Claim}
\newtheorem{remark}{Remark}
\nolinenumbers

\headers{Whittaker functions with one or both parameters large}{T. M. Dunster}

\title{Whittaker functions with one or both parameters large: simplified uniform asymptotic expansions involving Bessel and Airy functions}

\author{T. M. Dunster\thanks{Department of Mathematics and Statistics, San Diego State University, 5500 Campanile Drive, San Diego, CA 92182-7720, USA. 
  (\email{mdunster@sdsu.edu}, \url{https://tmdunster.sdsu.edu}).}
  }

\usepackage{amsopn}

\makeatletter
\newcommand*{\addFileDependency}[1]{
  \typeout{(#1)}
  \@addtofilelist{#1}
  \IfFileExists{#1}{}{\typeout{No file #1.}}
}
\makeatother

\nolinenumbers

\ifpdf
\hypersetup{
  pdftitle={Whittaker functions with one or both parameters large: simplified uniform asymptotic expansions involving Bessel and Airy functions},
  pdfauthor={T. M. Dunster}
}
\fi

\begin{document}

\maketitle

\begin{abstract}
Uniform asymptotic expansions are derived as $\kappa \to \infty$ for the Whittaker functions $W_{\kappa,\mu}(z)$, $M_{\kappa,\mu}(z)$, as well as related functions including generalized Laguerre polynomials. The results are uniformly valid for $0\leq\mu\leq(1-\delta_0)\kappa<\kappa$, where $\delta_0\in(0,1)$ is arbitrary and fixed. The analysis is based on the associated differential equation, which has a double pole and two turning points. At one of the turning points, which may coalesce with the double pole, a recently developed asymptotic theory is applied to obtain expansions involving Bessel functions. At the second turning point, expansions involving Airy functions are obtained. In both cases, the coefficients are readily computable, in contrast to those occurring in earlier results. The expansions, when taken together, uniformly cover the entire complex $z$-plane on the principal branch. Standard analytic continuation and connection formulas extend the results to all branches of $z$ and to negative $\mu$.
\end{abstract}



\section{Introduction}

We study the differential equation
\begin{equation}
\label{eq01}
\frac{d^2 w}{dt^2}
=\left\{
\frac{1}{4}
-\frac{\kappa}{t}
+\frac{\mu^2-\frac14}{t^2}
\right\}w,
\end{equation}
whose solutions include the Whittaker functions $\mathbf{M}_{\kappa,\mu}(t)$ and $W_{\kappa,\mu}(t)$. In terms of confluent hypergeometric functions \cite[Eqs.~13.14.2, 13.14.3]{NIST:DLMF}, these are given by
\begin{equation}
\label{eq02}
\mathbf{M}_{\kappa,\mu}(t)
=\frac{M_{\kappa,\mu}(t)}{\Gamma(1+2\mu)}
=\frac{t^{\mu+1/2}e^{-t/2}}
{\Gamma(1+2\mu)}
M\left(\mu-\kappa+\tfrac{1}{2};1+2\mu;t\right)
\end{equation}
and
\begin{equation}
\label{eq03}
W_{\kappa,\mu}(t)
=
t^{\mu+1/2}e^{-t/2}
U\left(\mu-\kappa+\tfrac12;1+2\mu;t\right).
\end{equation}
Here $\mathbf{M}_{\kappa,\mu}(t)$ denotes the Whittaker $M$-function scaled similarly to \cite[Eq.~13.2.4]{NIST:DLMF}; this normalization is convenient since the standard Whittaker $M$-function is not defined when $2\mu$ is a negative integer. Our purpose is to obtain uniform asymptotic expansions as $\kappa\to\infty$ that are valid for a large range of $\mu$ and throughout unbounded regions of the complex $z$-plane, with coefficients that are readily computable.

The classical Whittaker and confluent hypergeometric functions arise in numerous applications throughout mathematics and physics. Representative examples include exactly solvable Schr\"odinger equations and quantum bound-state problems \cite{Kubota:1969:SES,Nogueira:2016:LCB}, Coulomb wave functions and Green functions \cite{Esparza:1999:ZWF,Stovicek:2024:CGW}, singular quantum potentials \cite{Saad:2003:ICH}, wave propagation and optics \cite{Heading:1965:RIP}, and fluid stability \cite{Dyson:1960:SIA}; see also the classical monographs \cite{Buchholz:1969:CHF,Slater:1960:CHF}.

Large-parameter uniform asymptotic expansions for Whittaker and confluent hypergeometric functions have a long history, beginning with the pioneering work of Kazarinoff \cite{Kazarinoff:1955:AEW,Kazarinoff:1957:AFW} and Skovgaard \cite{Skovgaard:1966:UAE}, and continuing through later developments of Temme \cite{Temme:1978:UAC} and Dunster \cite{Dunster:2021:UAW}. Of particular relevance here are the results of Olver \cite{Olver:1980:WFW} and Dunster \cite{Dunster:1989:UAE}, which allow both Whittaker parameters to vary over substantial ranges. Olver obtained uniform approximations in terms of parabolic cylinder functions for $\mu\to\infty$, uniformly with respect to $x\in(0,\infty)$ and $\kappa\in[-(1-\delta)\mu,\mu]$, where $\delta>0$ is fixed; these results provide leading-order approximations, together with error bounds, rather than full asymptotic expansions. The present paper updates \cite{Dunster:1989:UAE} by deriving substantially simpler uniform expansions whose coefficients can be computed explicitly and efficiently, in contrast to the coefficients in the earlier expansions, which are generally not practicable to compute.

We first record the characteristic recessive behavior of the Whittaker functions that will be used later to identify the asymptotic solutions; these formulas follow, for example, from \cite[Eqs.~13.14.7, 13.14.14, 13.14.21, 13.19.2, 13.19.3]{NIST:DLMF}. As $t\to0$,
\begin{equation}
\label{eq04}
\mathbf{M}_{\kappa,\mu}(t)
=
\frac{t^{\mu+1/2}}{\Gamma(1+2\mu)}
\{1+\mathcal{O}(t)\}
\quad (t\to0,\, 2\mu \ne -1,-2,-3,\ldots),
\end{equation}
while in the exceptional cases
\begin{equation}
\label{eq05}
\mathbf{M}_{\kappa,-n/2}(t)
=
(-1)^n
\frac{\Gamma(\kappa+\frac12 n+\frac12)}
{n!\,\Gamma(\kappa-\frac12 n+\frac12)}
t^{(n+1)/2}\{1+\mathcal{O}(t)\}
\quad (t\to0,\, n = 1,2,3,\ldots).
\end{equation}
Thus, $\mathbf{M}_{\kappa,\mu}(t)$ is recessive at $t=0$ when $\Re(\mu) \ge 0$, since all other independent solutions are $\mathcal{O}(t^{-\mu+1/2})$ as $t \to 0$ ($\Re(\mu) > 0$), and $\mathcal{O}(t^{1/2}\ln(t))$ as $t \to 0$ ($\mu=0$).

Next, we have
\begin{equation}
\label{eq06}
W_{\kappa,\mu}(t)
=t^\kappa e^{-t/2}\left\{1+\mathcal{O}(t^{-1})\right\}
\quad
(t\to\infty,
\,
-\tfrac{3}{2}\pi<\arg(t)<\tfrac{3}{2}\pi),
\end{equation}
which is recessive as $t\to\infty$ in the right half-plane $|\arg(t)|<\tfrac12\pi$, since all other independent solutions grow exponentially there.

The following two analytic continuations provide important solutions that are recessive as $\Re(t) \to -\infty$ in the upper and lower half-planes of the principal plane:
\begin{equation}
\label{eq07}
W_{-\kappa,\mu}(te^{\pi i})
=e^{-\kappa\pi i}t^{-\kappa}e^{t/2}\left\{1+\mathcal{O}(t^{-1})\right\}
\quad
(t\to\infty,
\,
-\tfrac{5}{2}\pi<\arg(t)<\tfrac{1}{2}\pi)
\end{equation}
and
\begin{equation}
\label{eq08}
W_{-\kappa,\mu}(te^{-\pi i})
=e^{\kappa\pi i}t^{-\kappa}e^{t/2}\left\{1+\mathcal{O}(t^{-1})\right\}
\quad
(t\to\infty,
\,
-\tfrac{1}{2}\pi<\arg(t)<\tfrac{5}{2}\pi).
\end{equation}
We also note that
\begin{multline}
\label{eq09}
\mathbf{M}_{\kappa,\mu}(t)
=\frac{e^{(\mu-\kappa+1/2)\pi i}}
{\Gamma(\mu+\kappa+\tfrac{1}{2})}
t^\kappa e^{-t/2}\left\{1+\mathcal{O}(t^{-1})\right\}
\\
+
\frac{1}
{\Gamma(\mu-\kappa+\tfrac{1}{2})}
t^{-\kappa}e^{t/2}\left\{1+\mathcal{O}(t^{-1})\right\}
\quad
(t\to\infty,
\,
-\tfrac{1}{2}\pi<\arg(t)<\tfrac{3}{2}\pi).
\end{multline}

We shall also require the following connection formulas; see \cite[\S13.14(vii)]{NIST:DLMF}:
\begin{multline}
\label{eq10}
W_{\kappa,\mu}(t)
=\frac{1}{2\pi}
\Gamma\left(\kappa+\mu+\tfrac{1}{2}\right)
\Gamma\left(\kappa-\mu+\tfrac{1}{2}\right)
\\
\times
\left\{
e^{(\kappa-1/2)\pi i}
W_{-\kappa,\mu}(te^{\pi i})
+
e^{-(\kappa-1/2)\pi i}
W_{-\kappa,\mu}(te^{-\pi i})
\right\},
\end{multline}
\begin{equation}
\label{eq11}
\mathbf{M}_{\kappa,\mu}(t)
=\frac{1}{2\pi}
\Gamma\left(\kappa-\mu+\tfrac{1}{2}\right)
\left\{
e^{\mu\pi i}
W_{-\kappa,\mu}(te^{\pi i})
+
e^{-\mu\pi i}
W_{-\kappa,\mu}(te^{-\pi i})
\right\},
\end{equation}
\begin{equation}
\label{eq12}
\mathbf{M}_{\kappa,-\mu}(t)
=\frac{1}{2\pi}
\Gamma\left(\kappa+\mu+\tfrac{1}{2}\right)
\left\{
e^{-\mu\pi i}
W_{-\kappa,\mu}(te^{\pi i})
+
e^{\mu\pi i}
W_{-\kappa,\mu}(te^{-\pi i})
\right\}
\end{equation}
and
\begin{equation}
\label{eq13}
W_{\kappa,-\mu}(t)=W_{\kappa,\mu}(t).
\end{equation}
We note the analytic continuation formula \cite[Eq.~13.14.10]{NIST:DLMF}
\begin{equation}
\label{eq14}
\mathbf{M}_{-\kappa,\mu}(te^{\pm\pi i})
=e^{\pm(\mu+1/2)\pi i}\mathbf{M}_{\kappa,\mu}(t).
\end{equation}

From \cref{eq02} and \cite[Eqs.~13.14.26, 13.14.27 and 13.14.30]{NIST:DLMF}, the Wronskians are given by
\begin{equation}
\label{eq15}
\mathscr{W}\!\left\{\mathbf{M}_{\kappa,\mu}(t),
W_{\kappa,\mu}(t)\right\}
=
-\frac{1}
{\Gamma\!\left(\mu-\kappa+\tfrac12\right)},
\end{equation}
\begin{equation}
\label{eq16}
\mathscr{W}\!\left\{\mathbf{M}_{\kappa,\mu}(t),
W_{-\kappa,\mu}\!\left(t e^{\pm\pi i}\right)\right\}
=\mp i
\frac{e^{\mp\mu\pi i}}
{\Gamma\!\left(\mu+\kappa+\tfrac12\right)},
\end{equation}
\begin{equation}
\label{eq17}
\mathscr{W}\!\left\{W_{\kappa,\mu}(t),
W_{-\kappa,\mu}\!\left(t e^{\pm\pi i}
\right)\right\}
=
e^{\mp\kappa\pi i}
\end{equation}
and
\begin{equation}
\label{eq18}
\mathscr{W}\!\left\{
W_{-\kappa,\mu}\!\left(te^{\pi i}\right),
W_{-\kappa,\mu}\!\left(te^{-\pi i}\right)
\right\}
=
\frac{2\pi i}
{\Gamma\!\left(\kappa+\mu+\tfrac12\right)
\Gamma\!\left(\kappa-\mu+\tfrac12\right)}.
\end{equation}
Thus, for $0\le\mu<\kappa$, as assumed throughout this paper, the four solutions are pairwise linearly independent, with the exception of $\mathbf{M}_{\kappa,\mu}(t)$ and $W_{\kappa,\mu}(t)$ when $\kappa-\mu-\tfrac12=0,1,2,\ldots$.

We now introduce the appropriate scaling for our large $\kappa$ analysis. Let
\begin{equation}
\label{eq19}
z=t/\kappa, \quad \alpha=2\mu/\kappa,
\quad \sigma=(4-\alpha^2)^{1/2},
\end{equation}
so that \cref{eq01} takes the form
\begin{equation}
\label{eq20}
\frac{d^2 w}{dz^2}
=\left\{
\kappa^2 f(\alpha,z)
+g(z)
\right\}w,
\end{equation}
where
\begin{equation}
\label{eq21}
f(\alpha,z)
=\frac{z^2-4z+\alpha^2}{4z^2}
=\frac{(z-z_t)(z-\check{z}_t)}{4z^2},
\quad
g(z)=-\frac{1}{4z^2},
\end{equation}
with
\begin{equation}
\label{eq22}
z_t=2-\sigma,
\quad
\check{z}_t=2+\sigma.
\end{equation}
The point $z=0$ is a regular singularity of \cref{eq20}; indeed,
$f(\alpha,z)\sim\tfrac14\alpha^2z^{-2}$ as $z\to0$, and the corresponding exponents are $\tfrac12\pm\tfrac12\alpha\kappa=\tfrac12\pm\mu$. The zeros $z_t$ and $\check{z}_t$ of $f(\alpha,z)$ are the two turning points of the differential equation.

Throughout the paper we assume that, for some fixed $\sigma_0>0$,
\begin{equation}
\label{eq23}
0 < \sigma_0 \le \sigma \le 2,
\end{equation}
or, equivalently,
\begin{equation}
\label{eq24}
0 \le \alpha \le
\alpha_0 :=\sqrt{4-\sigma_0^2}<2.
\end{equation}
It follows from \cref{eq22,eq23} that
$0\le z_t\le2-\sigma_0<2+\sigma_0\le\check{z}_t\le4$,
so that the two turning points remain bounded away from each other. On the other hand, from \cref{eq19}, $\sigma\to2$ as $\alpha\to0$, in which case $z_t\to0$. Thus, the turning point $z_t$ is allowed to coalesce with the pole at $z=0$, an important feature of the present analysis. The lower bound in \cref{eq23} excludes the coalescence of the two turning points. The complementary limiting case $\sigma\to0$, equivalently $\alpha\to2$, requires parabolic cylinder functions as approximants and will be considered in a subsequent paper using the general theory developed in \cite{Dunster:2025:LCT}.

In terms of the original Whittaker parameters, \cref{eq19,eq24} show that the resulting asymptotic expansions hold as $\kappa\to\infty$, uniformly for the large parameter range
\begin{equation}
\label{eq25}
0\le\mu\le(1-\delta_0)\kappa<\kappa,
\end{equation}
where $\delta_0\in(0,1)$ is an arbitrary fixed constant.

The plan of the paper is as follows. In \cref{sec:Bessel} we apply the theory of \cite{Dunster:2026:TPD} to obtain Bessel-type expansions for solutions of \cref{eq20}, valid at both the pole $z=0$ and the turning point $z=z_t$, uniformly for the parameter range \cref{eq25}. In particular, the turning point is permitted to coalesce with the pole. These solutions are then matched with the standard Whittaker functions, with the principal results stated in \cref{thm:Bessel}.

The Bessel-type expansions are not valid at the second turning point $z=\check{z}_t$. In \cref{sec:Airy} we therefore obtain expansions that are valid there, again uniformly for the parameter range \cref{eq25}, using the theory developed in \cite{Dunster:2017:COA}. These involve Airy functions and their derivatives, as in the classical theory of \cite[Chap.~11]{Olver:1997:ASF}, but have substantially simpler and readily computable coefficients. The resulting solutions are matched with the standard Whittaker functions, and the main results are given in \cref{thm:Airy}.

Numerical examples are also given in \cref{sec:Bessel,sec:Airy} to illustrate the straightforward computation of the coefficients and the accuracy of the resulting expansions over their respective regions of validity. When taken together, the Bessel and Airy expansions are valid throughout the entire principal $z$-plane, and all results extend beyond it by means of the well-known analytic continuation formulas \cite[\S13.14(ii)]{NIST:DLMF}. 

Finally, in \cref{sec:Laguerre} we apply the results of the preceding sections to obtain uniform asymptotic expansions for the generalized Laguerre polynomials $L_n^{(a)}(z)$ as $n\to\infty$. The expansions are uniformly valid for
$-(1-\delta)n\le a\le\Delta n$,
where $\delta$ and $\Delta$ are arbitrary fixed positive constants. As for the Whittaker functions, the expansions involve Bessel and Airy functions and, when used together, are uniformly valid throughout the complex $z$-plane.

\section{Bessel expansions}
\label{sec:Bessel}

In this section we construct expansions that are valid in a domain containing both the pole $z=0$ and the turning point $z=z_t$, including the limiting case in which these points coalesce. We apply the theory developed in \cite{Dunster:2026:TPD}. It follows from \cref{eq21} that conditions (1.2) and (1.3) of that paper are satisfied, ensuring that the resulting approximations remain valid when the turning point coalesces with the pole. The construction also applies as $z\to\infty$, since the required conditions on $f(\alpha,z)$ and $g(z)$ are satisfied there; see \cite[Chap.~10, Exercise~4.1]{Olver:1997:ASF}.

The method involves the classical Liouville-Green (LG) variable $\xi$ (see, for example, \cite[Chap.~10]{Olver:1997:ASF}), together with a variable $\zeta$ that occurs in the arguments of the Bessel-function approximants. An important feature of the latter is that, unlike $\xi$, it is analytic at both $z=0$ and $z=z_t$. In the present case, \cite[Eqs.~(1.5), (1.10)-(1.12)]{Dunster:2026:TPD} give
\begin{equation}
\label{eq26}
\frac{d\zeta}{dz}
=\frac{\zeta}{z}
\left(\frac{\left(z_t-z\right)
\left(\check{z}_t-z\right)}
{\alpha^2-\zeta}\right)^{1/2}
\end{equation}
and
\begin{equation}
\label{eq27}
\int_{\alpha^2}^{\zeta}\frac{(\alpha^2-s)^{1/2}}
{2s}\,ds
=\xi
=\int_{z_t}^{z}
\frac{(z_t-v)^{1/2}(\check{z}_t-v)^{1/2}}{2v}\,dv.
\end{equation}
The branches are chosen so that $\zeta\in[0,\check{\zeta}_t]$ for $z\in[0,\check{z}_t]$ and $i\xi>0$ for $z\in(z_t,\check{z}_t)$, and are defined elsewhere by continuity with the appropriate cuts specified below.

Under this transformation, the pole $z=0$ corresponds to $\zeta=0$, while the turning point $z=z_t$ corresponds to $\zeta=\alpha^2$. Thus, as $\alpha\to0$ and $z_t\to0$, these two corresponding points in the $\zeta$-plane coalesce at the origin. The second turning point $z=\check{z}_t$ corresponds to $\zeta=\check{\zeta}_t$, where
\begin{equation}
\label{eq28}
\int_{\alpha^2}^{\check{\zeta}_t}
\frac{(s-\alpha^2)^{1/2}}
{2s}\,ds
=\int_{z_t}^{\check{z}_t}
\frac{(v-z_t)^{1/2}(\check{z}_t-v)^{1/2}}{2v}\,dv.
\end{equation}

We emphasize that $\zeta$ is analytic as a function of $z$ at both $z=0$ and $z=z_t$, but not at the second turning point $z=\check{z}_t$. By contrast, the LG variable $\xi$ is singular at $z=0$ and $z=z_t$.

\begin{figure}
 \centering
 \includegraphics[
 width=1.0\textwidth,keepaspectratio]{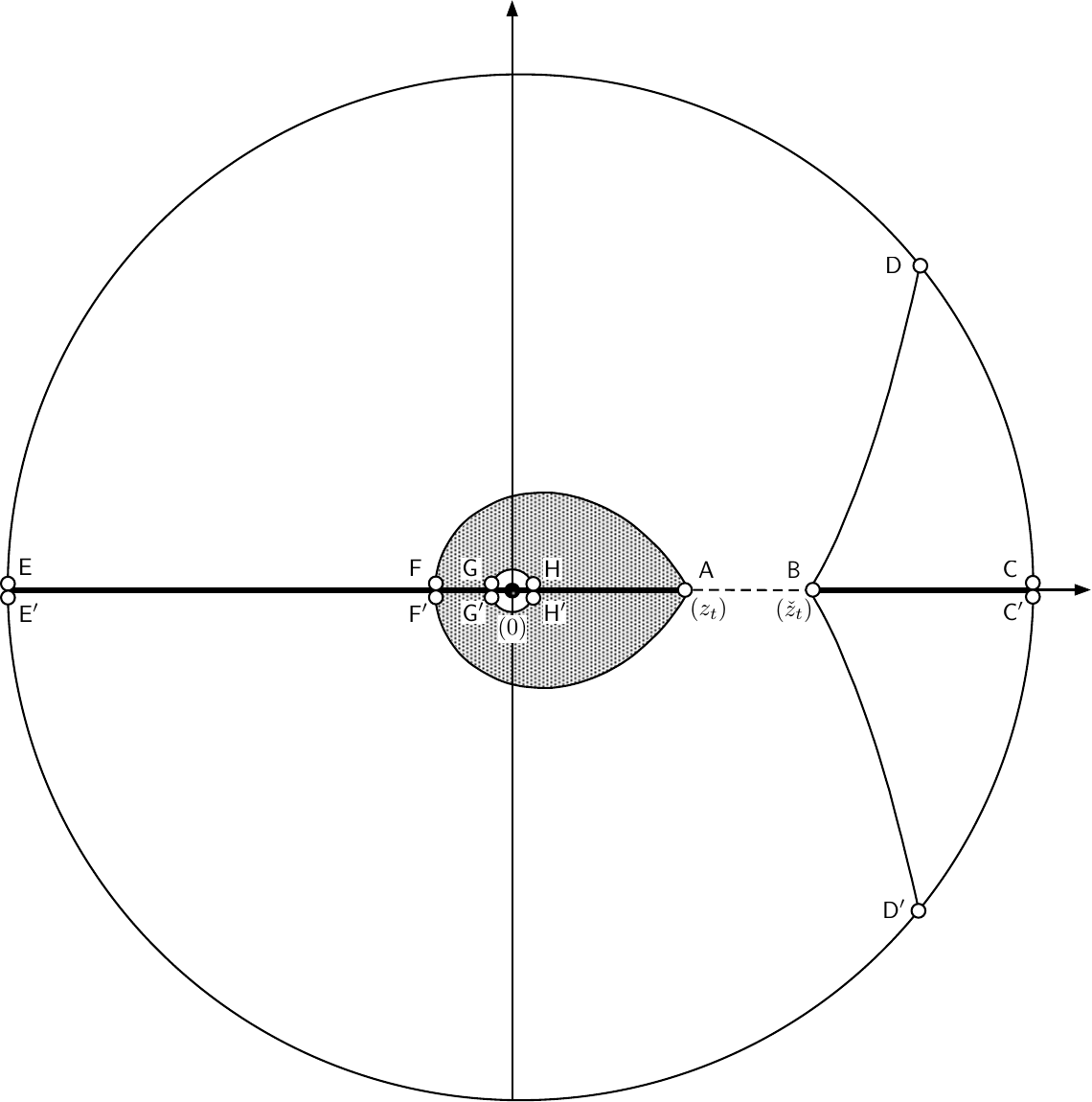}
 \caption{$z$-plane}
 \label{fig:zplane}
\end{figure}

\begin{figure}
 \centering
 \includegraphics[
 width=1.0\textwidth,keepaspectratio]{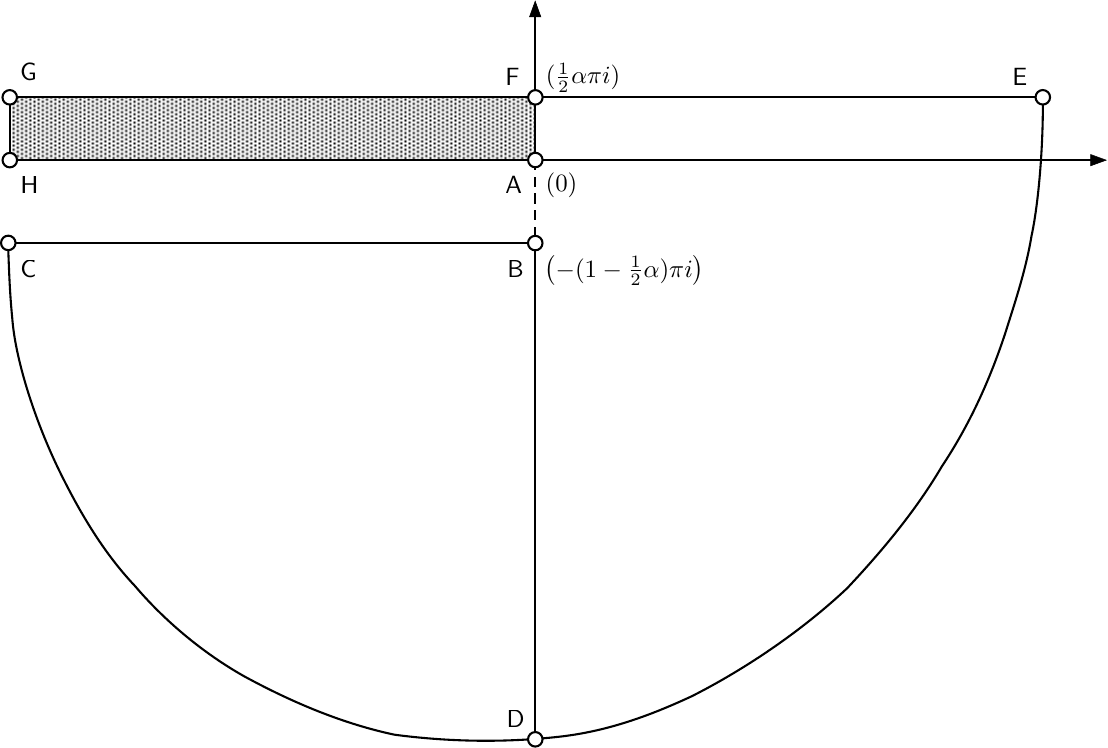}
 \caption{$\xi$ map of upper half $z$-plane}
 \label{fig:xiplane1}
\end{figure}

\begin{figure}
 \centering
 \includegraphics[
 width=1.0\textwidth,keepaspectratio]{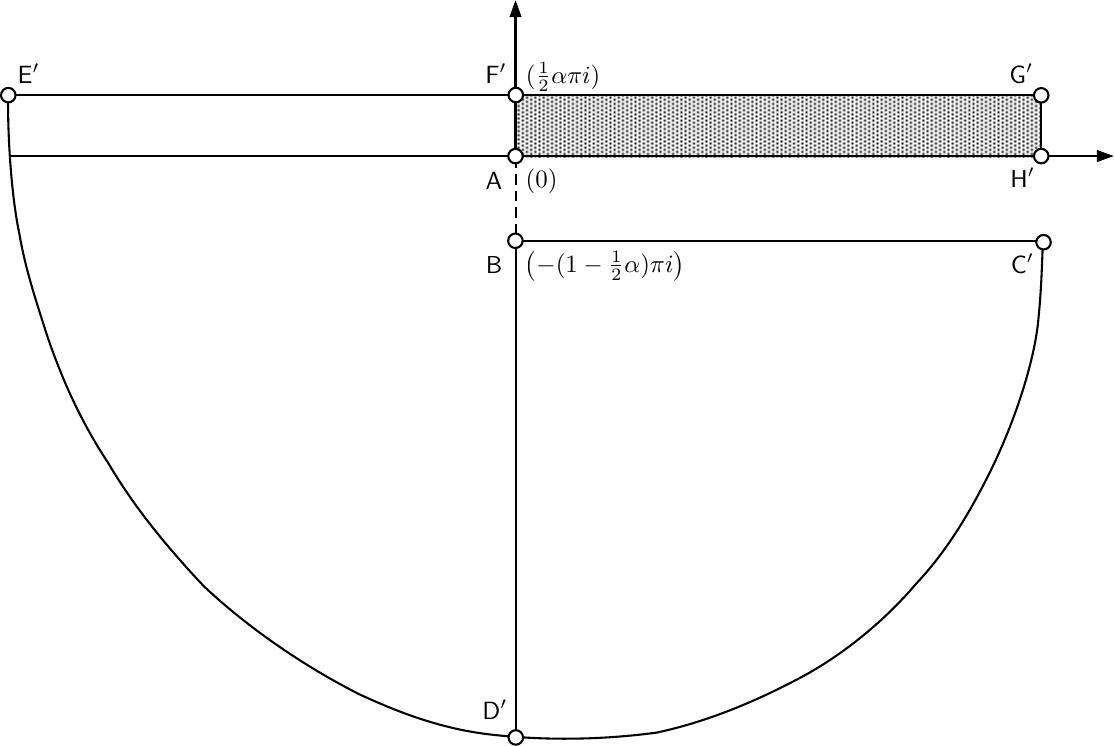}
 \caption{$\xi$ map of lower half $z$-plane}
 \label{fig:xiplane2}
\end{figure}

For later use we record explicit forms of the transformation. On integrating \cref{eq27}, for $\zeta>\alpha^2$ and $z_t<z<\check z_t$ we obtain
\begin{multline}
\label{eq29}
i\xi=\left(\zeta-\alpha^{2}\right)^{1/2}
-\alpha \arctan\left\{
\frac{\left(\zeta-\alpha^{2}\right)^{1/2}}{\alpha}
\right\}
\\
=\arccos\left(\frac{2-z}{\sigma}\right)
-\frac{\alpha}{2}
\arccos\left(\frac{\alpha^{2}-2z}{\sigma z}\right)
+\frac12\left\{\left(z-z_t\right)
\left(\check z_t-z\right)\right\}^{1/2},
\end{multline}
whereas for $\zeta\le\alpha^2$ and $z<z_t$ with $\alpha>0$,
\begin{multline}
\label{eq30}
(\alpha^2-\zeta)^{1/2}
-\frac{\alpha}{2}\ln\left\{
\frac{\alpha+(\alpha^2-\zeta)^{1/2}}
{\alpha-(\alpha^2-\zeta)^{1/2}}
\right\}
\\
=
\frac{1}{2}(z_t-z)^{1/2}(\check{z}_t-z)^{1/2}
+\frac{\alpha}{2}\ln\left\{
\left( \frac{z}{z_t} \right)
\frac{\alpha^2-2z_t}
{\alpha^2-2z+\alpha(z_t-z)^{1/2}(\check{z}_t-z)^{1/2}}
\right\}
\\
+\ln\left\{
\frac{2-z_t}
{2-z-(z_t-z)^{1/2}(\check{z}_t-z)^{1/2}}
\right\}.
\end{multline}
These relations extend to complex $z$, $\zeta$, and $\xi$ by analytic continuation. We take the $z$-plane to have cuts along $(-\infty,0]$ and $[\check{z}_t,\infty)$ for the mapping $z\mapsto\zeta$, and along $(-\infty,0]$ and $[z_t,\infty)$ for $z\mapsto\xi$. The corresponding $\xi$-maps of the upper and lower halves of the cut $z$-plane are shown in \cref{fig:zplane,fig:xiplane1,fig:xiplane2}, with corresponding points indicated; the figures are not drawn to scale.

We also require the behavior of $\xi$ at infinity. From \cref{eq29}, as $z\to\infty$ in the upper half-plane,
\begin{equation}
\label{eq31}
\xi=-\frac{1}{2}z+\ln(2z)
+1-\frac{1}{4}\alpha\ln\left(\frac{2+\alpha}
{2-\alpha}\right)
-\frac{1}{2}\ln(4-\alpha^2)
-\left(1-\frac{1}{2}\alpha\right)\pi i
+\mathcal{O}\left(\frac{1}{z}\right),
\end{equation}
and in the lower half-plane,
\begin{equation}
\label{eq32}
\xi
=
\frac{1}{2}z-\ln(2z)
-1+\frac{1}{4}\alpha\ln\left(\frac{2+\alpha}
{2-\alpha}\right)
+\frac{1}{2}\ln(4-\alpha^2)
-\left(1-\frac{1}{2}\alpha\right)\pi i
+\mathcal{O}\left(\frac{1}{z}\right).
\end{equation}

We now construct the coefficients required for the Bessel-type expansions. The LG Schwarzian, defined as in \cite[Eq.~(1.7)]{Dunster:2026:TPD}, is in the present case
\begin{multline}
\label{eq33}
\psi(\alpha,z)=\frac{4f(\alpha,z)
f''(\alpha,z)-5f'(\alpha,z)^2}
{16f(\alpha,z)^3}
+\frac{g(z)}{f(\alpha,z)}
\\
=
-\frac{z\left\{z^3+4(1-\alpha^2)z
+4\alpha^2\right\}}
{\left\{(z_t-z)(\check{z}_t-z)\right\}^3}.
\end{multline}
Although this is singular at the turning points, it satisfies $\psi(\alpha,z)=\mathcal{O}(z^{-2})$ as $z\to\infty$. We also note that, for $\alpha>0$, $\psi(\alpha,z)=\mathcal{O}(z)$ as $z\to0$, whereas it is singular there when $\alpha=0$.

Following \cite[Eqs.~(1.11), (1.15)]{Dunster:2020:LGE}, define
\begin{equation}
\label{eq34}
F_{1}(\alpha,z)=\frac12 \psi(\alpha,z),
\quad
F_{2}(\alpha,z)=-\frac{1}{4f^{1/2}(\alpha,z)}
\frac{d \psi(\alpha,z)}{dz}
\end{equation}
and
\begin{multline}
\label{eq35}
F_{s+1}(\alpha,z)=-\frac{1}{2f^{1/2}(\alpha,z)}
\frac{d F_{s}(\alpha,z)}{dz}
-\frac{1}{2}\sum_{j=1}^{s-1}F_{j}(\alpha,z)
F_{s-j}(\alpha,z)
\\
\quad (s=2,3,4,\ldots).
\end{multline}
The odd LG coefficients are then defined, in accordance with \cite[Eq.~(1.14)]{Dunster:2020:LGE}, by
\begin{equation}
\label{eq36}
E_{2s+1}(\alpha,z)=\int_{-\infty+i0}^{z}
F_{2s+1}(\alpha,v)f^{1/2}(\alpha,v)dv
+\lambda_{2s+1}(\alpha)
\quad (s=0,1,2,\ldots).
\end{equation}
The integration constants $\lambda_{2s+1}(\alpha)$ are chosen so that the resulting Bessel-type expansions remain valid at the turning point $z=z_t$. In the notation of \cite[Sec.~2]{Dunster:2026:TPD}, the two relevant singularities are $z_1=-\infty+i0$ and $z_2=-\infty-i0$, and the required values are
\begin{equation}
\label{eq37}
\lambda_{2s+1}(\alpha)=
E_{2s+1}(\alpha,-\infty+i0)
=-E_{2s+1}(\alpha,-\infty-i0)
\quad (s=0,1,2,\ldots).
\end{equation}
We shall see below that this condition determines each $\lambda_{2s+1}(\alpha)$ explicitly.

Even LG coefficients are determined by the formal expansion (see \cite[Eq. (2.9)]{Dunster:2017:COA})
\begin{equation}
\label{eq38}
\sum\limits_{s=1}^{\infty }
\frac{E_{2s}(\alpha,z) }{\kappa^{2s}}
\sim -\frac{1}{2}\ln \left\{ 1+\sum\limits_{s=0}^{\infty }
\frac{F_{2s+1}(\alpha,z) }{\kappa^{2s+2}}\right\}.
\end{equation}
From \cref{eq21,eq33,eq34,eq35}, using $\psi(\alpha,z)=\mathcal{O}(z^{-2})$ as $z\to\infty$, it follows by induction that
$F_{2s+1}(\alpha,z)\to0$ as $z\to-\infty\pm i0$. Accordingly,
\begin{equation}
\label{eq39}
E_{2s}(\alpha,-\infty+i0)
=E_{2s}(\alpha,-\infty-i0)
=0
\quad  (s=1,2,3,\ldots).
\end{equation}

To facilitate obtaining the coefficients in an explicitly computable form, we introduce, as in \cite[Eq.~(5.23)]{Dunster:2026:TPD}, the variable
\begin{equation}
\label{eq40}
\beta=\left(\frac{z-z_t}{z-\check{z}_t}\right)^{1/2},
\end{equation}
with the branch chosen consistently with that of $\xi$, so that $i\beta>0$ for $z\in(z_t,\check{z}_t)$. Elsewhere $\beta$ is defined by continuity in the plane cut along $(-\infty,z_t]$ and $[\check{z}_t,\infty)$. In particular, $\beta\to\pm1$ as $z\to\infty\pm i0$. Solving \cref{eq40} for $z$ gives
\begin{equation}
\label{eq41}
z=\frac{z_t-\beta^2 \check{z}_t}{1-\beta^2}.
\end{equation}

Next, differentiating \cref{eq40} with respect to $\xi$ yields
\begin{equation}
\label{eq42}
G(\alpha,\beta):=\frac{d\beta}{d\xi}
=
\frac{(1-\beta^2)^2\{\beta^2\check{z}_t-z_t\}}
{\beta^2(z_t-\check{z}_t)^2},
\end{equation}
or, equivalently, from \cref{eq19,eq22},
\begin{equation}
\label{eq43}
G(\alpha,\beta)
=
\frac{(1-\beta^2)^2
\{(2+\sigma)\beta^2-(2-\sigma)\}}
{4\sigma^2\beta^2}.
\end{equation}
Now, let $\mathrm{E}_s(\alpha,\beta)=E_s(\alpha,z)$, where $z$ and $\beta$ are related by \cref{eq41}. Since $d\xi=d\beta/G(\alpha,\beta)$, it follows from \cref{eq33,eq34,eq36,eq40,eq41,eq43} that the first odd LG coefficient is
\begin{multline}
\label{eq44}
\mathrm{E}_1(\alpha,\beta)
=
-\int_{1}^{\beta}
\frac{5(\sigma+2)p^6-(\sigma+6)p^4
+(\sigma-6)p^2+5\sigma-10}
{32\sigma^2p^4}\,dp
+\lambda_1(\alpha)
\\
=
-\frac{
5(\sigma+2)\beta^6
-3(\sigma+6)\beta^4
+3(\sigma-6)\beta^2
+5(2-\sigma)}
{96\sigma^2\beta^3}
-\frac{1}{6\sigma^2}+\lambda_1(\alpha).
\end{multline}
Given that $\beta\to\pm1$ as $z\to\infty\pm i0$, the choice \cref{eq37} becomes
\begin{equation}
\label{eq45}
\lambda_{2s+1}(\alpha)
=\mathrm{E}_{2s+1}(\alpha,1)
=-\mathrm{E}_{2s+1}(\alpha,-1)
\quad (s=0,1,2,\ldots).
\end{equation}
Thus the integration constants are chosen so that each $\mathrm{E}_{2s+1}(\alpha,\beta)$ is an odd Laurent polynomial in $\beta$. In particular, applying \cref{eq45} to \cref{eq44} gives
\begin{equation}
\label{eq46}
\lambda_1(\alpha)=\frac{1}{6\sigma^2},
\end{equation}
and consequently,
\begin{equation}
\label{eq47}
\mathrm{E}_1(\alpha,\beta)
=
-\frac{
5(\sigma+2)\beta^6
-3(\sigma+6)\beta^4
+3(\sigma-6)\beta^2
+5(2-\sigma)}
{96\sigma^2\beta^3}.
\end{equation}

As mentioned above, the even coefficients are normalized to vanish at $z=\infty$ in either direction, or, equivalently, at $\beta=\pm1$. From \cref{eq36,eq38,eq42}, they may be determined from equating like inverse powers of $\kappa$ in the asymptotic expansion
\begin{equation}
\label{eq48}
\sum\limits_{s=1}^{\infty }
\frac{\mathrm{E}_{2s}(\alpha,\beta) }{\kappa^{2s}}
\sim -\frac{1}{2}\ln \left\{ 1+G(\alpha,\beta)
\sum\limits_{s=1}^{\infty}
\frac{1}{\kappa^{2s}}
\pdv{\mathrm{E}_{2s-1}(\alpha,\beta)}{\beta}\right\}
\quad (\kappa \to \infty).
\end{equation}
For example, from \cref{eq43,eq47,eq48} we find that
\begin{multline}
\label{eq49}
\mathrm{E}_2(\alpha,\beta)
=
\frac{(\beta^2-1)^2}
{256\sigma^4\beta^6}
\left\{(\sigma+2)\beta^2+\sigma-2\right\}
\\
\times
\left\{
5(\sigma+2)\beta^6
-(\sigma+6)\beta^4
+(6-\sigma)\beta^2
+5(\sigma-2)
\right\}.
\end{multline}

All subsequent odd and even coefficients can be generated recursively from \cref{eq47,eq48} together with
\begin{multline}
\label{eq50}
\mathrm{E}_{2s+1}(\alpha,\beta)
=
-\frac{1}{2}G(\alpha,\beta)
\pdv{\mathrm{E}_{2s}(\alpha,\beta)}{\beta}
\\
-\frac{1}{4}
\left\{
\int_{\beta_{0}}^{\beta}
+\int_{-\beta_{0}}^{\beta}
\right\}
G(\alpha,p)
\sum_{j=1}^{2s-1}
\pdv{\mathrm{E}_j(\alpha,p)}{p}
\pdv{\mathrm{E}_{2s-j}(\alpha,p)}{p}\,dp
\quad (s=1,2,3,\ldots),
\end{multline}
which follows from \cref{eq36,eq42,eq45}. Here $\beta_{0}$ is an arbitrary nonzero constant and the integration paths avoid the pole at $p=0$. By induction, \cref{eq50} ensures that each odd coefficient is an odd function of $\beta$, as required.

A second family of coefficients is also required. Unlike the preceding family, these quantities are independent of the particular functions $f(\alpha,z)$ and $g(z)$ in the differential equation. They arise from the Bessel-function approximants and are given in \cite[Sec.~3]{Dunster:2026:TPD}. We record their construction here. First, let
\begin{equation}
\label{eq51}
\hat{\beta}=\left(\alpha^{2}-\zeta\right)^{1/2},
\end{equation}
where the branch of the square root is chosen so that $i\hat{\beta}$ is positive for $\alpha^{2}<\zeta<\infty$, and is defined by continuity in the $\zeta$-plane cut along $(-\infty,\alpha^{2}]$. Thus $\hat{\beta}$ is positive when real above the cut and negative when real below it.

We then define
\begin{equation}
\label{eq52}
e_{1}(\alpha,\hat{\beta})
=\frac{3\hat{\beta}^{2}-5\alpha^{2}}
{24\hat{\beta}^{3}},
\end{equation}
\begin{equation}
\label{eq53}
e_{2}(\alpha,\hat{\beta})
=\frac{\left(\hat{\beta}^{2}-\alpha^{2}\right)
\left(\hat{\beta}^{2}-5\alpha^{2}\right)}{16\hat{\beta}^{6}},
\end{equation}
\begin{equation}
\label{eq54}
\tilde{e}_{1}(\alpha,\hat{\beta})
=-\frac{9\hat{\beta}^{2}-7\alpha^{2}}{24\hat{\beta}^{3}}
\end{equation}
and
\begin{equation}
\label{eq55}
\tilde{e}_{2}(\alpha,\hat{\beta})
=-\frac{\left(\hat{\beta}^{2}-\alpha^{2}\right)
\left(3\hat{\beta}^{2}-7\alpha^{2}\right)}
{16\hat{\beta}^{6}},
\end{equation}
with subsequent coefficients generated by
\begin{multline}
\label{eq56}
e_{s+1}(\alpha,\hat{\beta})=
\frac{\alpha^{2}-\hat{\beta}^{2}}
{2\hat{\beta}^{2}}
\pdv{e_{s}(\alpha,\hat{\beta})}{\hat{\beta}}
\\
+\frac{1}{2}\int_{\infty}^{\hat{\beta}}
\frac{\alpha^{2}-p^{2}}{p^{2}}
\sum\limits_{j=1}^{s-1}
\pdv{e_{j}(\alpha,p)}{p}
\pdv{e_{s-j}(\alpha,p)}{p}\,dp
\quad (s=2,3,4,\ldots),
\end{multline}
and with $e$ replaced throughout by $\tilde e$ for $\tilde e_s(\alpha,\hat{\beta})$ ($s=3,4,5,\ldots$).

Combining the two families of coefficients as in \cite[Sec.~4]{Dunster:2026:TPD}, we define
\begin{equation}
\label{eq57}
\mathcal{E}_{s}(\alpha,z)
=\mathrm{E}_s(\alpha,\beta)
+(-1)^{s}e_{s}(\alpha,\hat{\beta})
\end{equation}
and
\begin{equation}
\label{eq58}
\tilde{\mathcal{E}}_{s}(\alpha,z)
=\mathrm{E}_s(\alpha,\beta)+(-1)^{s}
\tilde{e}_{s}(\alpha,\hat{\beta}).
\end{equation}

Next, from \cite[Thm.~4.2]{Dunster:2026:TPD}, two solutions of \cref{eq20}, recessive at $z=-\infty+i0$ and $z=-\infty-i0$, respectively, have the LG behavior
\begin{equation}
\label{eq59}
w_{\nu}^{(1,2)}(\kappa,\alpha,z) \sim f^{-1/4}(\alpha,z)
e^{\mp \kappa \xi \mp \pi i/4}
\quad (z \to -\infty \pm i0),
\end{equation}
where
\begin{equation}
\label{eq60}
\nu=\kappa\alpha=2\mu.
\end{equation}
Here $f^{-1/4}(\alpha,z)$ denotes the branch for which $\{-f(\alpha,z)\}^{-1/4}$ is positive for $z\in(z_t,\check z_t)$ and which is defined by continuity in the $z$-plane cut along $(-\infty,z_t]$ and $[\check z_t,\infty)$. To avoid a conflict with the standard Whittaker parameter $\mu$ used in this paper, we have replaced the parameter $\mu=u\alpha$ of \cite{Dunster:2026:TPD} by $\nu$ as defined in \cref{eq60}; the parameter $\alpha$ has the same meaning in both papers.

Again from \cite[Thm.~4.2]{Dunster:2026:TPD}, these solutions satisfy the connection relation
\begin{equation}
\label{eq61}
w_{\nu}^{(0)}(\kappa,\alpha,z)=w_{\nu}^{(1)}(\kappa,\alpha,z)
+\gamma_{\nu}(\kappa,\alpha)w_{\nu}^{(2)}(\kappa,\alpha,z),
\end{equation}
where $w_{\nu}^{(0)}(\kappa,\alpha,z)$ is a solution that is recessive at $z=0$, and, in general, for arbitrary $n>0$,
\begin{equation}
\label{eq62}
\gamma_{\nu}(\kappa,\alpha)=\gamma_{-\nu}(\kappa,\alpha)
=1+\mathcal{O}\left(\kappa^{-n}\right)
\quad (\kappa \to \infty).
\end{equation}

The corresponding Bessel-function expansions are then given by
\begin{multline}
\label{eq63}
\gamma_{\nu}^{(j)}(\kappa,\alpha)w_{\nu}^{(j)}(\kappa,\alpha,z)
=\left(\frac{\zeta-\alpha^{2}}{f(\alpha,z)}
\right)^{1/4}
\Biggl\{ \mathcal{C}_{\nu}^{(j)}\left(\kappa\zeta^{1/2}\right)A_{\nu}(\kappa,\alpha,z) \Biggr.
\\ \Biggl.
+\frac{\zeta}{\kappa}\frac{\partial
\mathcal{C}_{\nu}^{(j)}\left(\kappa\zeta^{1/2}\right)}{\partial \zeta}
B_{\nu}(\kappa,\alpha,z) \Biggr\},
\end{multline}
where $\gamma_{\nu}^{(0,1)}(\kappa,\alpha)=1$, $\gamma_{\nu}^{(2)}(\kappa,\alpha)=\gamma_{\nu}(\kappa,\alpha)$, $\mathcal{C}_{\nu}^{(0)}=2J_{\nu}$, and $\mathcal{C}_{\nu}^{(1,2)}=H_{\nu}^{(1,2)}$. The coefficient functions $A_{\nu}(\kappa,\alpha,z)$ and $B_{\nu}(\kappa,\alpha,z)$ are analytic for $z\in Z$ (described below) and, as $\kappa\to\infty$, possess the asymptotic expansions
\begin{multline} 
\label{eq64}
A_{\nu}(\kappa,\alpha,z) \sim 
\sqrt{\frac{\kappa \pi}{2}}\exp \left\{\sum\limits_{s=0}^{\infty}
\frac{\lambda_{2s+1}(\alpha)}{\kappa^{2s+1}}\right\}
\\ \times
\exp \left\{ \sum\limits_{s=1}^{\infty}
\frac{\tilde{\mathcal{E}}_{2s}(\alpha,z) }{\kappa^{2s}}\right\} 
\cosh \left\{ \sum\limits_{s=0}^{\infty}
\frac{\tilde{\mathcal{E}}_{2s+1}(\alpha,z)}
{\kappa^{2s+1}}\right\}
\end{multline}
and
\begin{multline} 
\label{eq65}
B_{\nu}(\kappa,\alpha,z) 
\sim \frac{\sqrt{2\pi \kappa}}{\hat{\beta}}
\exp \left\{\sum\limits_{s=0}^{\infty}
\frac{\lambda_{2s+1}(\alpha)}{\kappa^{2s+1}}\right\}
\\ \times
\exp \left\{ \sum\limits_{s=1}^{\infty}
\frac{\mathcal{E}_{2s}(\alpha,z) }{\kappa^{2s}}\right\} 
\sinh \left\{ \sum\limits_{s=0}^{\infty}
\frac{\mathcal{E}_{2s+1}(\alpha,z)}
{\kappa^{2s+1}}\right\},
\end{multline}
uniformly for $z\in Z$ and $\alpha\in[0,\alpha_0]$.

The expansions \cref{eq64,eq65} are constructed from LG expansions associated with the solutions recessive at $z=-\infty+i0$ and $z=-\infty-i0$, together with the corresponding expansions of the Bessel-function approximants and their derivatives. Their region of validity is therefore the intersection of the LG regions associated with these two recessive solutions. The general form of such domains is described in \cite[Chap.~10]{Olver:1997:ASF}. Roughly speaking, they consist of points that can be joined to the appropriate recessive endpoints by suitable piecewise-smooth paths along which the relevant LG conditions are satisfied. Such paths must remain bounded away from the turning points, where the LG expansions cease to be valid.

We now describe this region more precisely in the present case. Let $\mathscr{P}_{-1,1}$ denote a path that joins $z=-\infty+i0$ to $z=-\infty-i0$ that

\begin{enumerate}
\item[(i)] consists of a finite chain of $R_2$-arcs, as defined in \cite[Chap.~5, \S3.4]{Olver:1997:ASF};
\item[(ii)] the continuous extension of $\Re(\xi)$ is monotonic as $z$ traverses the path from one endpoint to the other; and
\item[(iii)] does not pass through the turning points $z=z_t$ and $z=\check{z}_t$.
\end{enumerate}

The region of validity, denoted by $Z$, is the set of points that lie on at least one such path. From \cref{fig:zplane,fig:xiplane1,fig:xiplane2}, it is seen that the region lying on and to the right of the curve $\mathsf{DBD'}$ is excluded. Indeed, no path joining $z=-\infty+i0$ to $z=-\infty-i0$ through a point of this region can preserve the monotonicity of $\Re(\xi)$ without passing through the turning point $z=\check{z}_t$. The closed pear-shaped shaded region containing $z=0$ and bounded by $\mathsf{AFF'}$, which we denote by $Z_0$, is excluded for a similar reason: a path through this region satisfying the monotonicity condition would have to pass through the other turning point $z=z_t$; see \cref{fig:xiplane1,fig:xiplane2}.

In summary, the region of validity $Z$ consists of the points in the cut $z$-plane that exclude those lying on or to the right of the curve $\mathsf{DBD'}$ and those in $Z_0$.

The excluded region $Z_0$, which contains the pole $z=0$ and the turning point $z=z_t$, is of particular interest. One method of extending the expansions into $Z_0$ is to use Cauchy's integral formula; see \cite{Dunster:2017:COA} for details. Here we instead use an alternative method, re-expanding \cref{eq64,eq65} in the conventional form as asymptotic series in inverse powers of $\kappa$. This is achieved by following the procedure described in \cite[Remark~1]{Dunster:2025:SAR}. Thus, let $\mathrm{A}_{s}(\alpha,z)=\tilde{q}_{2s}(\alpha,z)$ ($s=1,2,3,\ldots$) and $\mathrm{B}_{s}(\alpha,z)=\hat{\beta}^{-1}q_{2s+1}(\alpha,z)$ ($s=0,1,2,\ldots$), where, for $s=1,2,3,\ldots$,
\begin{equation}
\label{eq66}
\tilde{q}_{s}(\alpha,z)=\tilde{\mathcal{E}}_{s}(\alpha,z)
+\frac{1}{s}\sum_{j=1}^{s-1}j
\tilde{\mathcal{E}}_{j}(\alpha,z)\tilde{q}_{s-j}(\alpha,z)
\end{equation}
and
\begin{equation}
\label{eq67}
q_{s}(\alpha,z)=\mathcal{E}_{s}(\alpha,z)
+\frac{1}{s}\sum_{j=1}^{s-1}j
\mathcal{E}_{j}(\alpha,z)q_{s-j}(\alpha,z),
\end{equation}
with both sums empty when $s=1$. It then follows from \cite[Thm.~4.2]{Dunster:2026:TPD} that the following expansions hold throughout $Z\cup Z_0$, in particular at the turning point and the pole, uniformly for $\alpha\in[0,\alpha_0]$ as $\kappa\to\infty$:
\begin{equation} 
\label{eq68}
A_{\nu}(\kappa,\alpha,z) \sim 
\sqrt{\frac{\kappa \pi}{2}}\exp \left\{\sum\limits_{s=0}^{\infty}
\frac{\lambda_{2s+1}(\alpha)}{\kappa^{2s+1}}\right\}
\left\{1 + \sum_{s=1}^{\infty}
\frac{\mathrm{A}_{s}(\alpha,z)}{\kappa^{2s}}\right\}
\end{equation}
and
\begin{equation} 
\label{eq69}
B_{\nu}(\kappa,\alpha,z) \sim 
\sqrt{\frac{2\pi}{\kappa}} 
\exp \left\{\sum\limits_{s=0}^{\infty}
\frac{\lambda_{2s+1}(\alpha)}{\kappa^{2s+1}}\right\}
\sum_{s=0}^{\infty}
\frac{\mathrm{B}_{s}(\alpha,z)}{\kappa^{2s}}.
\end{equation}

For $\alpha\in[0,\alpha_0]$, the coefficients $\mathrm{A}_{s}(\alpha,z)$ and $\mathrm{B}_{s}(\alpha,z)$ are analytic at $z=0$ and $z=z_t$, after removal of the apparent singularities there. From \cref{eq57,eq58,eq66,eq67}, we determine the leading ones as 
$\mathrm{A}_1(\alpha,z)=\tilde{\mathcal{E}}_2(\alpha,z)+\tfrac12\{\tilde{\mathcal{E}}_1(\alpha,z)\}^2$
and $\mathrm{B}_0(\alpha,z)=\hat{\beta}^{-1}\mathcal{E}_1(\alpha,z)$.

We now match these asymptotic solutions with the Whittaker functions. Since both sides are recessive solutions of \cref{eq20} as $z\to-\infty+i0$, a comparison of \cref{eq07,eq19,eq21,eq31,eq59} gives
\begin{equation}
\label{eq70}
W_{-\kappa,\mu}(\kappa z e^{-\pi i})
=\frac{e^{(\mu-1/4)\pi i}}{\sqrt{2}}
\left(\frac{2e}{\kappa(4-\alpha^2)^{1/2}}\right)^\kappa
\left(\frac{2-\alpha}{2+\alpha}\right)^{\mu/2}
w_{\nu}^{(1)}(\kappa,\alpha,z).
\end{equation}
Similarly, a comparison of \cref{eq08,eq19,eq21,eq32,eq59} for solutions that are recessive as $z\to-\infty-i0$ gives
\begin{equation}
\label{eq71}
W_{-\kappa,\mu}(\kappa z e^{\pi i})
=\frac{e^{-(\mu+1/4)\pi i}}{\sqrt{2}}
\left(\frac{2e}{\kappa(4-\alpha^2)^{1/2}}\right)^\kappa
\left(\frac{2-\alpha}{2+\alpha}\right)^{\mu/2}
w_{\nu}^{(2)}(\kappa,\alpha,z).
\end{equation}
Then it follows from \cref{eq11,eq61,eq70,eq71} that
$\gamma_{\nu}(\kappa,\alpha)=1$. Hence, using \cref{eq10,eq11,eq60,eq63,eq70,eq71}, we obtain
\begin{multline}
\label{eq72}
W_{-\kappa,\mu}(\kappa z e^{\mp \pi i})
=e^{\pm \mu\pi i}
\left(\frac{2e}{\kappa(4-\alpha^2)^{1/2}}\right)^\kappa
\left(\frac{2-\alpha}{2+\alpha}\right)^{\mu/2}
\\  \times
z^{1/2}\left(\frac{\alpha^{2}-\zeta}
{(z_t-z)(\check{z}_t-z)}\right)^{1/4}
\Biggl\{ H_{2\mu}^{(1,2)}\left(\kappa\zeta^{1/2}\right)A_{\nu}(\kappa,\alpha,z) 
\\
+\frac{\zeta}{\kappa}\frac{\partial
H_{2\mu}^{(1,2)}\left(\kappa\zeta^{1/2}\right)}{\partial \zeta}
B_{\nu}(\kappa,\alpha,z) \Biggr\},
\end{multline}
where the branch of the quarter-power factor is positive for $z\in(0,\check{z}_t)$.

\subsection{Main Bessel expansions and numerical results}

We now present the principal results of this section. Let $\alpha$, $z_t$, $\check{z}_t$, $\zeta$, and $\xi$ be given by \cref{eq19,eq22,eq27}, and let $\beta$ and $\hat{\beta}$ be defined by \cref{eq40,eq51}. The coefficients $\mathrm{E}_s(\alpha,\beta)$ and the constants $\lambda_{2s+1}(\alpha)$ are determined by \cref{eq45,eq47,eq48,eq50}, while $e_s(\alpha,\hat{\beta})$ and $\tilde{e}_s(\alpha,\hat{\beta})$ are given by \cref{eq52,eq53,eq54,eq55,eq56}. The combinations $\mathcal{E}_s(\alpha,z)$ and $\tilde{\mathcal{E}}_s(\alpha,z)$ are then defined by \cref{eq57,eq58}. Finally, let $\tilde{q}_s(\alpha,z)$ and $q_s(\alpha,z)$ be given by \cref{eq66,eq67}, and hence $\mathrm{A}_s(\alpha,z)$ and $\mathrm{B}_s(\alpha,z)$ as specified immediately before \cref{eq66}. The formulas below are obtained from \cref{eq10}, \cref{eq11}, \cref{eq12,eq72}.

\begin{theorem}
\label{thm:Bessel}
\begin{multline}
\label{eq73}
\mathbf{M}_{\kappa,\mu}(\kappa z)
=\frac{\Gamma\left(\kappa-\mu+\tfrac{1}{2}\right)}
{\pi}
\left(\frac{2e}{\kappa(4-\alpha^2)^{1/2}}\right)^\kappa
\left(\frac{2-\alpha}{2+\alpha}\right)^{\mu/2} z^{1/2}
\\  \times
\left(\frac{\alpha^{2}-\zeta}
{(z_t-z)(\check{z}_t-z)}\right)^{1/4}
\Biggl\{ J_{2\mu}\left(\kappa\zeta^{1/2}\right)
A_{\nu}(\kappa,\alpha,z) 
+\frac{\zeta}{\kappa}\frac{\partial
J_{2\mu}\left(\kappa\zeta^{1/2}\right)}{\partial \zeta}
B_{\nu}(\kappa,\alpha,z) \Biggr\},
\end{multline}
\begin{multline}
\label{eq74}
\mathbf{M}_{\kappa,-\mu}(\kappa z)
=\frac{\Gamma\left(\kappa+\mu+\tfrac{1}{2}\right)}
{\pi}
\left(\frac{2e}{\kappa(4-\alpha^2)^{1/2}}\right)^\kappa
\left(\frac{2-\alpha}{2+\alpha}\right)^{\mu/2} z^{1/2}
\\ \times
\left(\frac{\alpha^{2}-\zeta}
{(z_t-z)(\check{z}_t-z)}\right)^{1/4}
\Biggl\{
\mathscr{C}_{2\mu}\left(\kappa\zeta^{1/2},2\mu\pi\right)
A_{\nu}(\kappa,\alpha,z)
\\
+\frac{\zeta}{\kappa}
\frac{\partial}{\partial\zeta}
\left\{
\mathscr{C}_{2\mu}\left(\kappa\zeta^{1/2},2\mu\pi\right)
\right\}
B_{\nu}(\kappa,\alpha,z)
\Biggr\}
\end{multline}
and
\begin{multline}
\label{eq75}
W_{\kappa,\mu}(\kappa z)
=\frac{\Gamma(\kappa+\mu+\tfrac{1}{2})
\Gamma\left(\kappa-\mu+\tfrac{1}{2}\right)}
{\pi}
\left(\frac{2e}{\kappa(4-\alpha^2)^{1/2}}\right)^\kappa
\left(\frac{2-\alpha}{2+\alpha}\right)^{\mu/2}
\\  \times
z^{1/2}
\left(\frac{\alpha^{2}-\zeta}
{(z_t-z)(\check{z}_t-z)}\right)^{1/4}
\Biggl\{ \mathscr{C}_{2\mu}\left(\kappa\zeta^{1/2},
\left(\mu-\kappa+\tfrac12\right)\pi\right)
A_{\nu}(\kappa,\alpha,z) 
\Biggr.
\\
\Biggl.
+\frac{\zeta}{\kappa}\frac{\partial}
{\partial \zeta}
\left\{
\mathscr{C}_{2\mu}\left(\kappa\zeta^{1/2},
\left(\mu-\kappa+\tfrac12\right)\pi\right)
\right\}
B_{\nu}(\kappa,\alpha,z) \Biggr\},
\end{multline}
where
\begin{equation}
\label{eq76}
\mathscr{C}_{\nu}(x,\theta)
=
\cos(\theta)J_{\nu}(x)-\sin(\theta)Y_{\nu}(x).
\end{equation}
The functions $A_{\nu}(\kappa,\alpha,z)$ and $B_{\nu}(\kappa,\alpha,z)$ are analytic at $z=0$ and $z=z_t$, and possess the asymptotic expansions \cref{eq68,eq69} as $\kappa\to\infty$, uniformly for $0\le\mu\le(1-\delta_0)\kappa<\kappa$, $\Re(z)\le\check{z}_t-\delta_1$, and $-\pi\le\arg(z)\le\pi$, where $\delta_0 \in (0,1)$ and $\delta_1 \in (0,\check{z}_t-z_t)$ are arbitrary fixed constants.
\end{theorem}

\begin{remark}
The region of validity in the $z$-plane is larger than the one stated in \cref{thm:Bessel}, and is as described earlier in this section. For simplicity, in the theorem we restrict attention to a convenient half-plane subregion that contains $z=0$ and $z=z_t$.
\end{remark}

We now illustrate the accuracy of the Bessel expansion for $\mathbf{M}_{\kappa,\mu}(\kappa z)$, the solution recessive at $z=0$, for real $z\in(0,\check{z}_t)$. For this purpose, define the approximation $\mathsf{M}_4(\kappa,\mu,z)$ by retaining the first four terms of the series in \cref{eq68,eq69,eq73}, viz.
\begin{multline}
\label{eq77}
\mathsf{M}_4(\kappa,\mu,z)
=
\frac{\Gamma\left(\kappa-\mu+\tfrac{1}{2}\right)}
{\pi}
\left(\frac{2e}{\kappa(4-\alpha^2)^{1/2}}\right)^\kappa
\left(\frac{2-\alpha}{2+\alpha}\right)^{\mu/2}
\\  \times
z^{1/2}
\left(\frac{\alpha^{2}-\zeta}
{(z_t-z)(\check{z}_t-z)}\right)^{1/4}
\left(\frac{\kappa\pi}{2}\right)^{1/2}
\exp\left\{\sum\limits_{s=0}^{3}
\frac{\lambda_{2s+1}(\alpha)}{\kappa^{2s+1}}\right\}
\\ \times
\Biggl[
J_{2\mu}\left(\kappa\zeta^{1/2}\right)
\left\{
1+\sum\limits_{s=1}^{3}
\frac{\mathrm{A}_{s}(\alpha,z)}{\kappa^{2s}}
\right\}
+
\frac{2\zeta}{\kappa^{2}}
\frac{\partial
J_{2\mu}\left(\kappa\zeta^{1/2}\right)}
{\partial \zeta}
\sum\limits_{s=0}^{3}
\frac{\mathrm{B}_{s}(\alpha,z)}{\kappa^{2s}}
\Biggr].
\end{multline}
For numerical computation near the turning point $z=z_t$, the removable singularities in the coefficients $\mathrm{A}_{s}(\alpha,z)$ and $\mathrm{B}_{s}(\alpha,z)$ are conveniently handled by expanding them in powers of $\beta$ about $\beta=0$, using a procedure similar to that described in \cite[Sec.~5]{Dunster:2026:TPD}.

A straightforward relative error obtained by dividing the absolute error $|\mathbf{M}_{\kappa,\mu}(\kappa z)\allowbreak-\mathsf{M}_4(\kappa,\mu,z)|$ by $|\mathbf{M}_{\kappa,\mu}(\kappa z)|$ is not suitable, since $\mathbf{M}_{\kappa,\mu}(\kappa z)$ has zeros in the oscillatory part of the interval. We therefore introduce a continuous positive envelope function that is close to the amplitude of $\mathbf{M}_{\kappa,\mu}(\kappa z)$ in the oscillatory interval and coincides with it outside that interval.

To this end, in view of \cref{eq11}, we first define the numerically satisfactory solution
\begin{equation}
\label{eq78}
\mathbf{N}_{\kappa,\mu}(t)
=\frac{\Gamma\left(\kappa-\mu+\tfrac12\right)}
{2\pi i}
\left\{e^{-\mu\pi i}W_{-\kappa,\mu}(t e^{-\pi i})
-e^{\mu\pi i}W_{-\kappa,\mu}(t e^{\pi i})
\right\}.
\end{equation}
It follows from \cref{eq72} that
\begin{multline}
\label{eq79}
\mathbf{N}_{\kappa,\mu}(\kappa z)
=
\frac{\Gamma\left(\kappa-\mu+\tfrac{1}{2}\right)}
{\pi}
\left(\frac{2e}{\kappa(4-\alpha^2)^{1/2}}\right)^\kappa
\left(\frac{2-\alpha}{2+\alpha}\right)^{\mu/2}
z^{1/2}
\\  \times
\left(\frac{\alpha^{2}-\zeta}
{(z_t-z)(\check{z}_t-z)}\right)^{1/4}
\left\{
Y_{2\mu}\left(\kappa\zeta^{1/2}\right)
A_{\nu}(\kappa,\alpha,z)
+\frac{\zeta}{\kappa}
\frac{\partial
Y_{2\mu}\left(\kappa\zeta^{1/2}\right)}
{\partial \zeta}
B_{\nu}(\kappa,\alpha,z)
\right\},
\end{multline}
and comparison of this with \cref{eq73} verifies that $\mathbf{N}_{\kappa,\mu}(\kappa z)$ is the appropriate companion. We then choose our envelope with the desired properties to be
\begin{equation}
\label{eq80}
\mathsf{M}(\kappa,\mu,z)
=
\begin{cases}
\mathbf{M}_{\kappa,\mu}(\kappa z)
& (0<z<\omega_{\kappa,\mu}),\\[5pt]
\left[
\left\{\mathbf{M}_{\kappa,\mu}(\kappa z)\right\}^2
+
\left\{\mathbf{N}_{\kappa,\mu}(\kappa z)\right\}^2
\right]^{1/2}
& (\omega_{\kappa,\mu}\le z<\check z_t),
\end{cases}
\end{equation}
where $\omega_{\kappa,\mu}$ is the smallest positive zero, in the $z$-variable, of $\mathbf{N}_{\kappa,\mu}(\kappa z)$. On $(0,\omega_{\kappa,\mu})$ the function $\mathbf{M}_{\kappa,\mu}(\kappa z)$ is positive and nonoscillatory, while in $[\omega_{\kappa,\mu},\check z_t)$ the second expression provides a positive envelope that closely tracks its oscillatory amplitude.

\begin{figure}[!htb]
 \centering
 \includegraphics[
 width=0.9\textwidth,keepaspectratio]{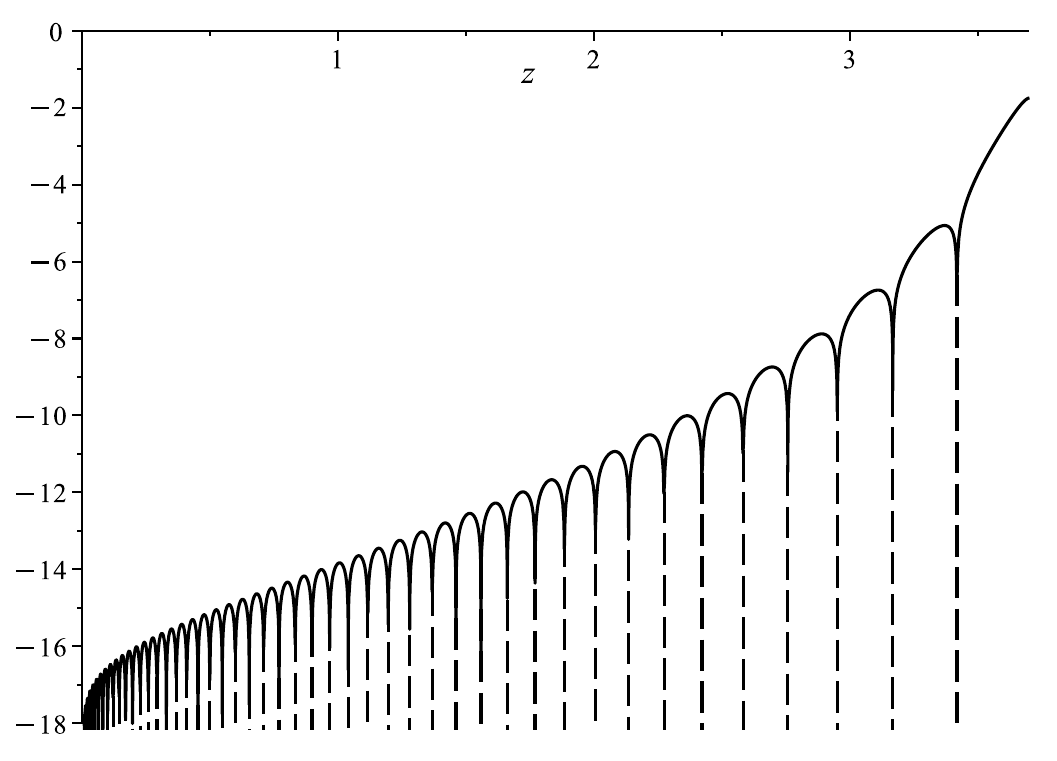}
 \caption{Graph of $\Omega_{4}(\kappa,\mu,z)$ for $\kappa=50$, $\alpha=0.01$ ($\mu=0.25$)}
 \label{fig:omega0.01}
\end{figure}

\begin{figure}[!htb]
 \centering
 \includegraphics[
 width=0.9\textwidth,keepaspectratio]{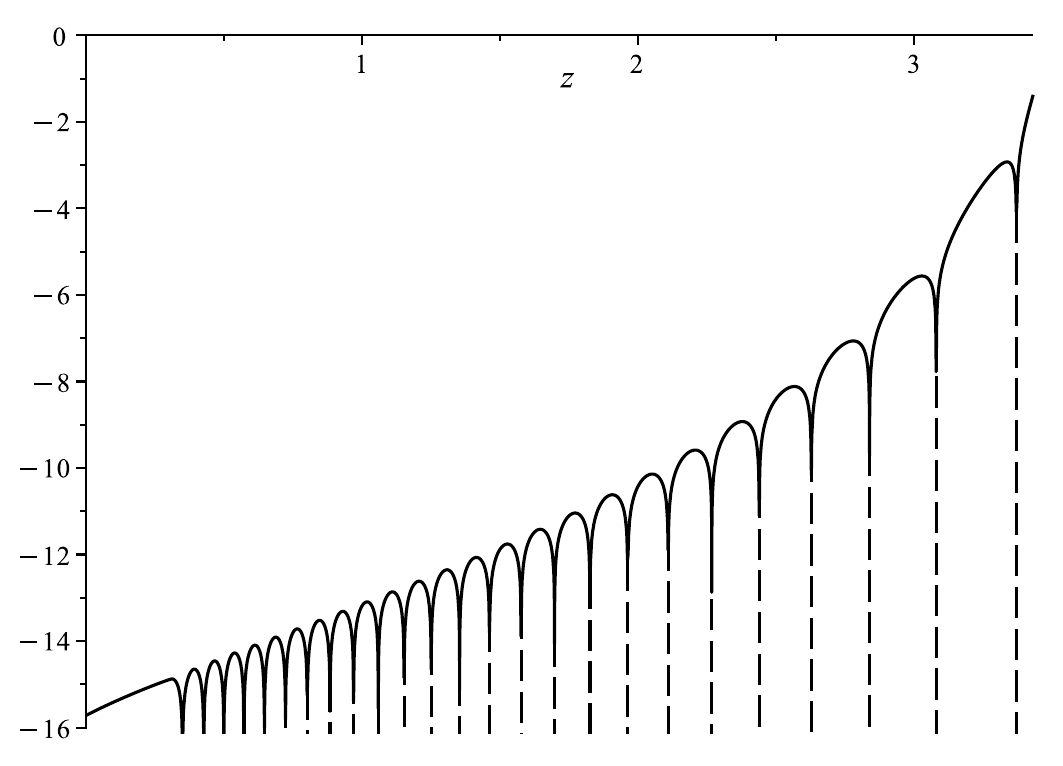}
 \caption{Graph of $\Omega_{4}(\kappa,\mu,z)$ for $\kappa=50$, $\alpha=1$ ($\mu=25$)}
 \label{fig:omega1}
\end{figure}

\begin{figure}[!htb]
 \centering
 \includegraphics[
 width=0.9\textwidth,keepaspectratio]{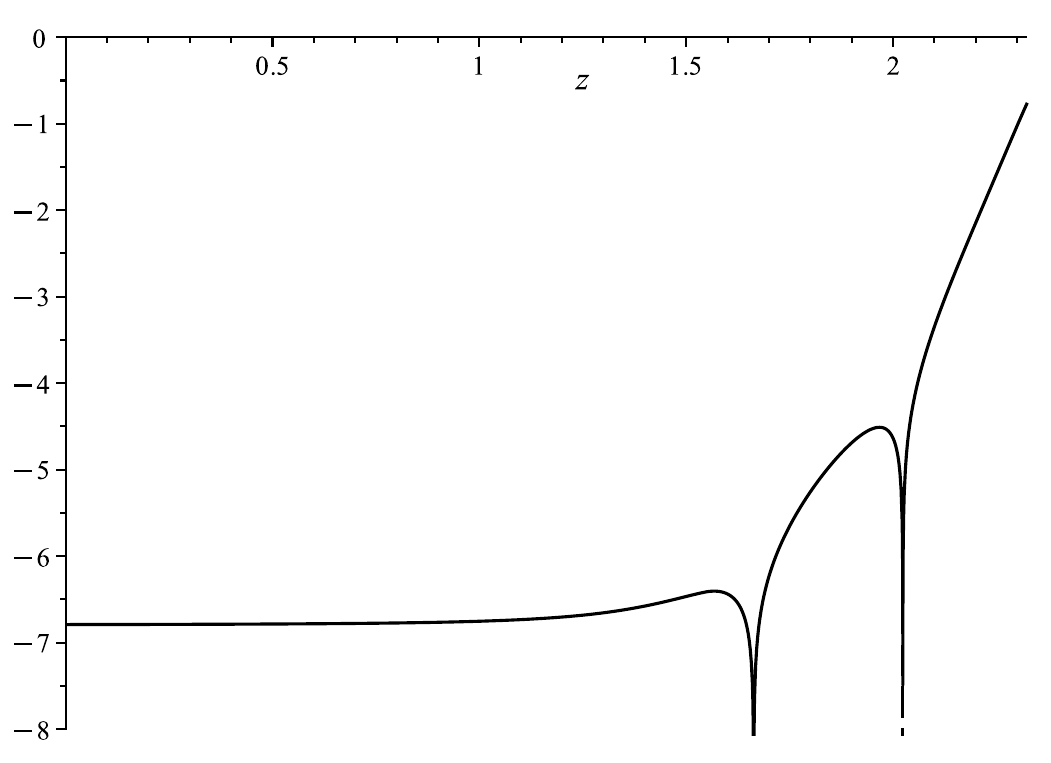}
\caption{Graph of $\Omega_{4}(\kappa,\mu,z)$ for $\kappa=50$, $\alpha=1.9$ ($\mu=47.5$)}
 \label{fig:omega1.9}
\end{figure}

The relative error was then measured by
\begin{equation}
\label{eq81}
\Omega_{4}(\kappa,\mu,z)
=
\log_{10}\left\{
\frac{
\left|\mathbf{M}_{\kappa,\mu}(\kappa z)
-\mathsf{M}_4(\kappa,\mu,z)\right|
}
{\mathsf{M}(\kappa,\mu,z)}
\right\}
\quad (0<z<\check{z}_t).
\end{equation}
For $\kappa=50$, plots of $\Omega_4(\kappa,\mu,z)$ are shown in
\cref{fig:omega0.01,fig:omega1,fig:omega1.9} for
$\alpha=0.01$, $1$, and $1.9$, respectively. The corresponding turning points are approximately
$z_t=2.5\times10^{-5}$, $0.27$, and $1.4$, and
$\check z_t=4.0$, $3.7$, and $2.6$, respectively. The sharp downward cusps in the graphs occur near points in the oscillatory interval where the approximation error $\mathbf{M}_{\kappa,\mu}(\kappa z)-\mathsf{M}_4(\kappa,\mu,z)$ vanishes.


In all three cases the accuracy deteriorates as $z$ approaches the second turning point $\check z_t$, as expected, since the Bessel expansion is not uniform there. For $\alpha=0.01$, the turning point $z_t$ lies extremely close to the pole at $z=0$, yet the approximation remains highly accurate in their neighborhood (16 digits or better) and performs well uniformly away from $\check z_t$. Comparable precision is obtained for $\alpha=1$, when $z_t$ is bounded away both from the pole and the second turning point $\check z_t$. When $\alpha=1.9$, however, the two turning points are much closer together, and a corresponding loss of precision is evident throughout the interval, although the approximation remains useful except near $\check z_t$. In this case an expansion based on coalescing turning points and parabolic cylinder functions would be more appropriate; as mentioned earlier, this case will be treated in a subsequent paper.

\section{Airy expansions}
\label{sec:Airy}

In this section we construct expansions that are valid at the second turning point $z=\check{z}_t$. The classical theory for a simple turning point gives expansions in terms of Airy functions; see \cite[Chap.~11]{Olver:1997:ASF}. Here we use the more recent theory developed in \cite{Dunster:2017:COA}, which likewise involves Airy functions and their derivatives, but has the advantage that the coefficients are readily computable.

The corresponding Airy variable $\check{\zeta}$ and associated LG variable $\check{\xi}$ are defined in the standard manner by
\begin{multline}
\label{eq82}
\check{\xi}
=\frac{2}{3}\check{\zeta}^{3/2}
= \int_{\check{z}_t}^{z} f^{1/2}(\alpha,t)\,dt
\\
=
\frac{1}{2}
\left(z-z_t\right)^{1/2}
\left(z-\check{z}_t\right)^{1/2}
-\ln\left\{
\frac{
z-2+
\left(z-z_t\right)^{1/2}
\left(z-\check{z}_t\right)^{1/2}
}
{\sigma}
\right\}
\\
+
\frac{\alpha}{4}
\ln\left\{
\frac{2z-\alpha^2
+\alpha
\left(z-z_t\right)^{1/2}
\left(z-\check{z}_t\right)^{1/2}
}
{2z-\alpha^2
-\alpha
\left(z-z_t\right)^{1/2}
\left(z-\check{z}_t\right)^{1/2}
}
\right\}.
\end{multline}

For $z\in(z_t,\check{z}_t]$, where $\check{\zeta}\le0$, an equivalent real form is
\begin{multline}
\label{eq83}
\frac{2}{3}(-\check{\zeta})^{3/2}
=
\int_{z}^{\check{z}_t}\{-f(\alpha,t)\}^{1/2}\,dt
\\
=
\arccos\left(\frac{z-2}{\sigma}\right)
-\frac12\left(z-z_t\right)^{1/2}
\left(\check{z}_t-z\right)^{1/2}
-\alpha
\arctan\left\{
\left(
\frac{(2-\sigma)(\check z_t-z)}
{(2+\sigma)(z-z_t)}
\right)^{1/2}
\right\}.
\end{multline}
The branches are chosen so that $\check{\zeta}>0$ for $z\in(\check{z}_t,\infty)$, and $\check{\zeta}$ is defined elsewhere by continuity in the $z$-plane cut along $(-\infty,z_t]$. It is analytic at the turning point $z=\check{z}_t$, where $\check{\zeta}=0$. The corresponding branch of $\check{\xi}=\frac{2}{3}\check{\zeta}^{3/2}$ is positive for $z\in(\check{z}_t,\infty)$ and is defined by continuity in the $z$-plane cut along $(-\infty,\check{z}_t]$.

\begin{figure}
 \centering
 \includegraphics[
 width=1.0\textwidth,keepaspectratio]{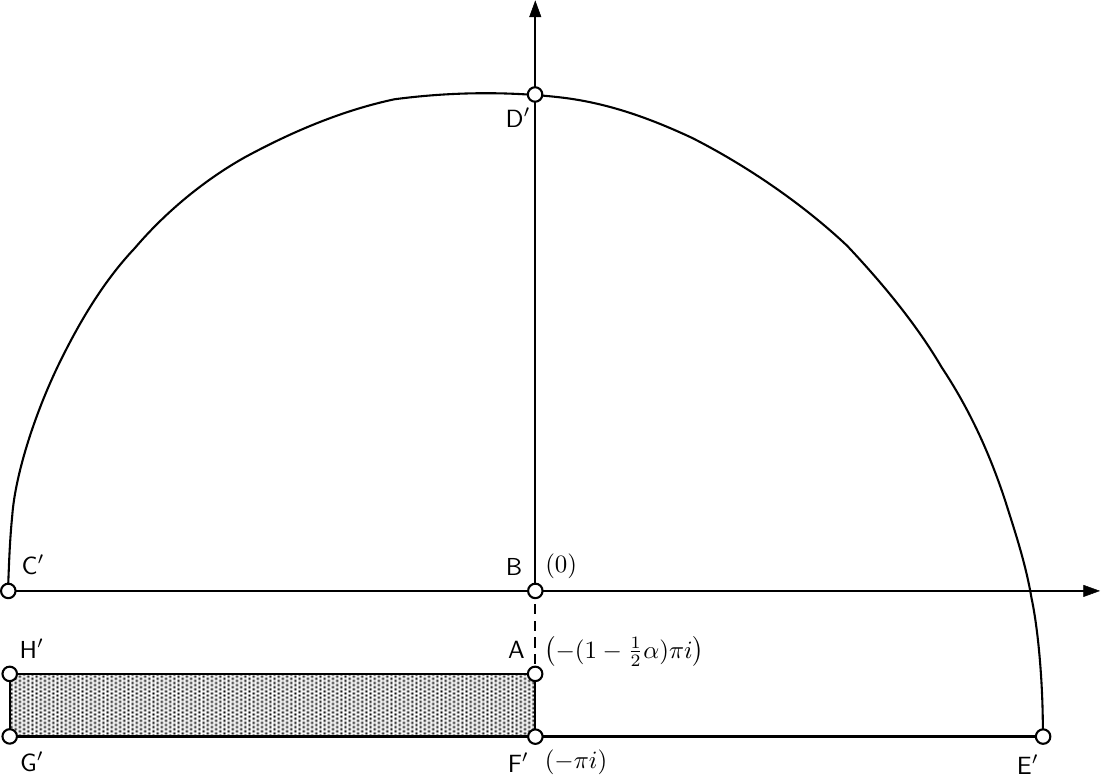}
 \caption{$\check{\xi}_{\mathsf{AB}}$ map of lower half $z$-plane}
 \label{fig:hatxiplane2}
\end{figure}

For the description of the region of validity below, we also require the branch $\check{\xi}_{\mathsf{AB}}$ obtained by continuation from the upper half $z$-plane across the interval $[z_t,\check{z}_t]$, labeled $\mathsf{AB}$ in \cref{fig:zplane}. This branch is related to $\check{\xi}$ by
\begin{equation}
\label{eq84}
\check{\xi}_{\mathsf{AB}}
=
\begin{cases}
\check{\xi} & (\Im(z)\ge 0),\\
-\check{\xi} & (\Im(z)<0).
\end{cases}
\end{equation}
The map of $\check{\xi}_{\mathsf{AB}}$ for the lower half $z$-plane is shown in \cref{fig:hatxiplane2}. From \cref{eq27,eq82,eq84},
\begin{equation}
\label{eq85}
\check{\xi}_{\mathsf{AB}}
=-\xi-\left(1-\tfrac12\alpha\right)\pi i.
\end{equation}

We shall also require the behavior at infinity, and from \cref{eq82},
\begin{equation}
\label{eq86}
\check{\xi}
=\frac{1}{2}z-\ln(2z)
-1+\frac{1}{4}\alpha\ln\left(\frac{2+\alpha}
{2-\alpha}\right)
+\frac{1}{2}\ln(4-\alpha^2)
+\mathcal{O}\left(\frac{1}{z}\right)
\quad (z\to\infty).
\end{equation}
\begin{figure}
 \centering
 \includegraphics[
 width=1.0\textwidth,keepaspectratio]{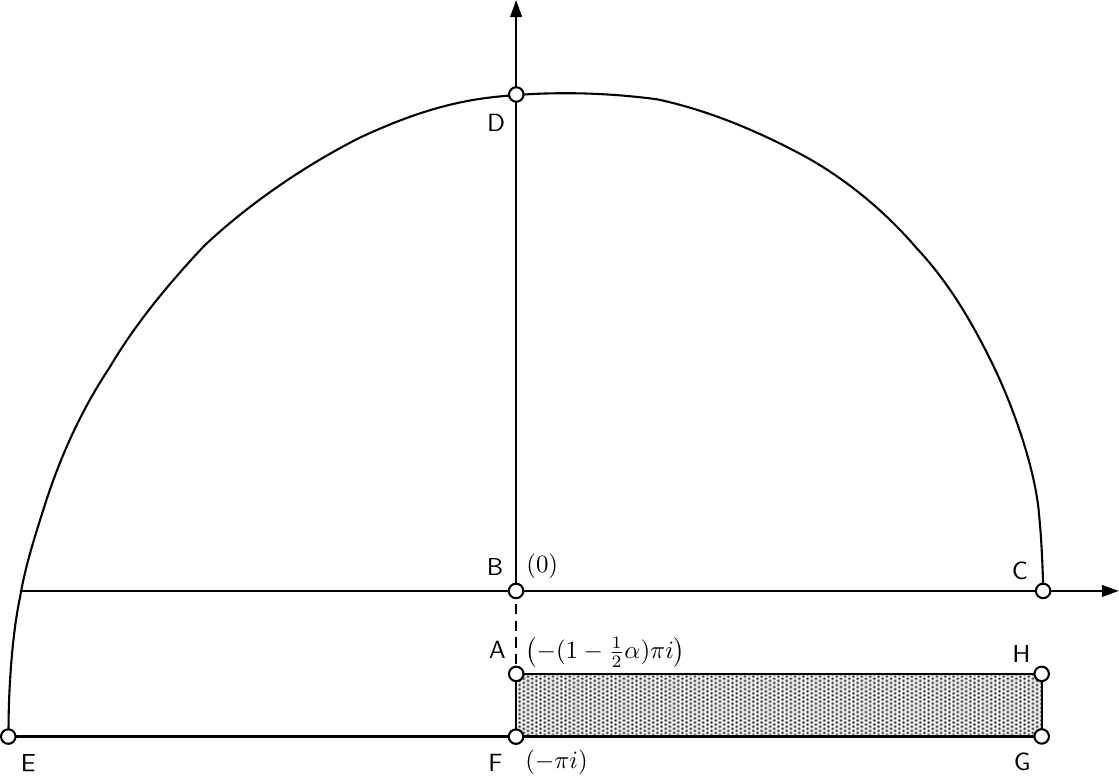}
 \caption{$\check{\xi}$ map of upper half $z$-plane}
 \label{fig:hatxiplane1}
\end{figure}

The three LG solutions required in the Airy construction are characterized by the asymptotic behaviors\footnote{Here and below, some of the notation differs from that of \cite{Dunster:2017:COA,Dunster:2021:SEB} in order to conform with that used elsewhere in this paper.}
\begin{equation}
\label{eq87}
\check{w}_{0}(\kappa,\alpha,\check{\zeta})
\sim
f^{-1/4}(\alpha,z)
e^{-\kappa\check{\xi}}
\quad (\check{\xi} \to +\infty)
\end{equation}
and
\begin{equation}
\label{eq88}
\check{w}_{\pm 1}(\kappa,\alpha,\check{\zeta})
\sim
f^{-1/4}(\alpha,z)
e^{\kappa\check{\xi}}
\quad (\check{\xi} \to -\infty \pm i0),
\end{equation}
where $f^{-1/4}(\alpha,z)$ denotes the branch that is positive for $z>\check{z}_t$ and is defined by continuity in the $z$-plane cut along $(-\infty,\check{z}_t]$. These limits uniquely determine the solutions by their recessive behavior at the respective singularities.

From \cite[Eqs.~(2.2), (2.3), (2.9), (2.13)]{Dunster:2021:SEB}, the three LG solutions satisfy the connection relation
\begin{equation}
\label{eq89}
\check{w}_{0}(\kappa,\alpha,\check{\zeta})=
i\check{c}_1(\kappa,\alpha)\check{w}_{1}(\kappa,\alpha,\check{\zeta})
-i\check{c}_{-1}(\kappa,\alpha)\check{w}_{-1}(\kappa,\alpha,\check{\zeta}),
\end{equation}
where the connection coefficients satisfy
\begin{equation}
\label{eq90}
\check{c}_{\pm 1}(\kappa,\alpha)
\exp\left\{
-\sum\limits_{s=1}^{\infty}
\frac{\lambda_s^{(\pm 1)}(\alpha)}{\kappa^s}
\right\}
\sim
\exp\left\{
-\sum\limits_{s=1}^{\infty}
(-1)^s\frac{\lambda_s^{(0)}(\alpha)}{\kappa^s}
\right\}
\quad (\kappa \to \infty),
\end{equation}
with
\begin{equation}
\label{eq91}
\lambda_s^{(0)}(\alpha)
=\lim_{\check{\xi} \to \infty} E_{s}(\alpha,z)
\end{equation}
and
\begin{equation}
\label{eq92}
\lambda_s^{(\pm 1)}(\alpha)
=\lim_{\check{\xi} \to -\infty \pm i 0} E_{s}(\alpha,z).
\end{equation}
The LG coefficients $E_s(\alpha,z)$ are generated from \cref{eq33,eq34,eq35}, with
\begin{equation}
\label{eq93}
E_{s}(\alpha,z)
=\int F_{s}(\alpha,z)
f^{1/2}(\alpha,z)\,dz,
\end{equation}
where the constants of integration are chosen so that $(z-\check{z}_t)^{1/2}E_{2s+1}(\alpha,z)$ is meromorphic at $z=\check{z}_t$ for $s=0,1,2,\ldots$. These coefficients are similar, but not identical, to those used in \cref{sec:Bessel}.

We next introduce two recursively-defined sequences $\{a_s\}_{s=1}^{\infty}$ and $\{\tilde a_s\}_{s=1}^{\infty}$. They are initialized by $a_1=a_2=5/72$ and $\tilde a_1=\tilde a_2=-7/72$, with subsequent terms in either sequence determined by
\begin{equation}
\label{eq94}
b_{s+1}=\frac{1}{2}(s+1)b_s+\frac{1}{2}
\sum\limits_{j=1}^{s-1}b_jb_{s-j}
\quad (b=a,\tilde a).
\end{equation}
Define
\begin{equation}
\label{eq95}
\mathscr{E}_{s}(\alpha,z)
=
E_s(\alpha,z)
+(-1)^s a_s s^{-1}\check{\xi}^{-s}
\end{equation}
and
\begin{equation}
\label{eq96}
\tilde{\mathscr{E}}_{s}(\alpha,z)
=
E_s(\alpha,z)
+(-1)^s \tilde a_s s^{-1}\check{\xi}^{-s}.
\end{equation}
Then the Airy expansions of \cite{Dunster:2017:COA} give
\begin{multline}
\label{eq97}
\frac{1}{2\pi^{1/2}\kappa^{1/6}}
\check{w}_{0}(\kappa,\alpha,\check{\zeta})
=
\left(\frac{\check{\zeta}}
{f(\alpha,z)}\right)^{1/4}
\left\{
\mathrm{Ai}\left(\kappa^{2/3}\check{\zeta}\right)
\check{A}(\kappa,\alpha,z)
\right.
\\
\left.
+
\mathrm{Ai}'\left(\kappa^{2/3}\check{\zeta}\right)
\check{B}(\kappa,\alpha,z)\right\},
\end{multline}
\begin{multline}
\label{eq98}
\frac{e^{\pi i/6}\check{c}_1(\kappa,\alpha)}{2\pi^{1/2}\kappa^{1/6}}
\check{w}_{1}(\kappa,\alpha,\check{\zeta})
=
\left(\frac{\check{\zeta}}
{f(\alpha,z)}\right)^{1/4}
\left\{
\mathrm{Ai}_1\left(\kappa^{2/3}\check{\zeta}\right)
\check{A}(\kappa,\alpha,z)
\right.
\\
\left.
+
\mathrm{Ai}'_1\left(\kappa^{2/3}\check{\zeta}\right)
\check{B}(\kappa,\alpha,z)\right\}
\end{multline}
and
\begin{multline}
\label{eq99}
\frac{e^{-\pi i/6}\check{c}_{-1}(\kappa,\alpha)}{2\pi^{1/2}\kappa^{1/6}}
\check{w}_{-1}(\kappa,\alpha,\check{\zeta})
=
\left(\frac{\check{\zeta}}
{f(\alpha,z)}\right)^{1/4}
\left\{
\mathrm{Ai}_{-1}\left(\kappa^{2/3}\check{\zeta}\right)
\check{A}(\kappa,\alpha,z)
\right.
\\
\left.
+
\mathrm{Ai}'_{-1}\left(\kappa^{2/3}\check{\zeta}\right)
\check{B}(\kappa,\alpha,z)\right\},
\end{multline}
where $\mathrm{Ai}_{j}(z):=\mathrm{Ai}(ze^{-2\pi ij/3})$ ($j=0,\pm1$). The coefficient functions have the asymptotic expansions
\begin{multline}
\label{eq100}
\check{A}(\kappa,\alpha,z)
\sim
\exp\left\{
-\sum\limits_{s=1}^{\infty}
(-1)^s\frac{\lambda_s^{(0)}(\alpha)}{\kappa^s}
\right\}
\exp\left\{\sum\limits_{s=1}^{\infty}
\frac{\tilde{\mathscr{E}}_{2s}(\alpha,z)}{\kappa^{2s}}
\right\}
\\ \times
\cosh\left\{\sum\limits_{s=0}^{\infty}
\frac{\tilde{\mathscr{E}}_{2s+1}(\alpha,z)}{\kappa^{2s+1}}
\right\}
\end{multline}
and
\begin{multline}
\label{eq101}
\check{B}(\kappa,\alpha,z)
\sim
\frac{1}
{\kappa^{1/3}\check{\zeta}^{1/2}}
\exp\left\{
-\sum\limits_{s=1}^{\infty}
(-1)^s\frac{\lambda_s^{(0)}(\alpha)}{\kappa^s}
\right\}
\exp\left\{\sum\limits_{s=1}^{\infty}
\frac{\mathscr{E}_{2s}(\alpha,z)}{\kappa^{2s}}
\right\}
\\ \times
\sinh\left\{\sum\limits_{s=0}^{\infty}
\frac{\mathscr{E}_{2s+1}(\alpha,z)}{\kappa^{2s+1}}
\right\},
\end{multline}
as $\kappa\to\infty$. The region of validity is described using three types of paths, since the construction of \cite{Dunster:2017:COA} involves LG expansions associated with the three fundamental solutions $\check{w}_{j}(\kappa,\alpha,\check{\zeta})$ ($j=0,\pm1$), in contrast to the two required in \cref{sec:Bessel}.

For two of these three types, let $\check{\mathscr{P}}_{0,1}$ and $\check{\mathscr{P}}_{0,-1}$ denote paths joining $z=+\infty$ to $z=-\infty+i0$ and $z=-\infty-i0$, respectively, and satisfying conditions (i)-(iii) stated in \cref{sec:Bessel}, except that $\xi$ in condition (ii) is replaced by $\check{\xi}$. The third type, denoted by $\check{\mathscr{P}}_{-1,1}$, joins $z=-\infty+i0$ to $z=-\infty-i0$, with $\xi$ in condition (ii) now replaced by the branch $\check{\xi}_{\mathsf{AB}}$. The region of validity, denoted by $\check{Z}$, is the set of points in the $z$-plane lying on at least one such path.

We now describe $\check{Z}$ explicitly. Consider first the upper half of the cut $z$-plane. From \cref{fig:zplane,fig:hatxiplane1}, the monotonicity condition on $\Re(\check{\xi})$, together with the requirement that the paths remain bounded away from the turning points, shows that every point in this half-plane lies on a suitable path $\check{\mathscr{P}}_{0,1}$, except for those on the interval $[z_t,\check{z}_t]$ and those in the part of $Z_0$ lying in the upper half-plane. By conjugation, the corresponding result holds in the lower half $z$-plane for paths $\check{\mathscr{P}}_{0,-1}$. Thus these two types of paths cover every point in the cut $z$-plane except those on $[z_t,\check{z}_t]$ and those in $Z_0$.

Furthermore, every point of $(z_t,\check{z}_t)$ lies on a suitable path $\check{\mathscr{P}}_{-1,1}$. As shown in \cref{fig:hatxiplane1,fig:hatxiplane2}, such a path may be chosen to cross $\mathsf{AB}$ while preserving the monotonicity of $\Re(\check{\xi}_{\mathsf{AB}})$. Indeed, $\Re(\check{\xi}_{\mathsf{AB}})\to-\infty$ as $z\to-\infty+i0$, whereas $\Re(\check{\xi}_{\mathsf{AB}})\to+\infty$ as $z\to-\infty-i0$.

In summary, the region of validity $\check{Z}$ consists of all points in the cut $z$-plane except those in $Z_0$ and the turning point $z=\check{z}_t$. In contrast with the region $Z$ described in \cref{sec:Bessel}, points lying on and to the right of the curve $\mathsf{DBD'}$ are included here, since they lie on suitable paths $\check{\mathscr{P}}_{0,1}$ or $\check{\mathscr{P}}_{0,-1}$. We shall subsequently extend the expansions to include a full neighborhood of $z=\check{z}_t$.

Before making this extension, we match the asymptotic solutions with the Whittaker functions. From \cref{eq06,eq19,eq21,eq22,eq86,eq87}, comparison of the solutions recessive at $z=+\infty$ ($\check{\xi}=+\infty$) gives
\begin{equation}
\label{eq102}
W_{\kappa,\mu}(\kappa z)
=
\frac{1}{\sqrt{2}}
\left(\frac{\kappa(4-\alpha^2)^{1/2}}{2e}\right)^\kappa
\left(\frac{2+\alpha}{2-\alpha}\right)^{\mu/2}
\check{w}_{0}(\kappa,\alpha,\check{\zeta}).
\end{equation}

Similarly, from \cref{eq07,eq08,eq19,eq21,eq22,eq86,eq88}, comparison of the solutions recessive at $z=-\infty\pm i0$ gives
\begin{equation}
\label{eq103}
W_{-\kappa,\mu}(\kappa z e^{\mp \pi i})
=
\frac{e^{\pm\kappa\pi i}}{\sqrt{2}}
\left(\frac{2e}{\kappa(4-\alpha^2)^{1/2}}\right)^\kappa
\left(\frac{2-\alpha}{2+\alpha}\right)^{\mu/2}
\check{w}_{\pm 1}(\kappa,\alpha,\check{\zeta}).
\end{equation}
Substitution of \cref{eq102,eq103} into \cref{eq10}, followed by comparison with \cref{eq89}, yields
\begin{equation}
\label{eq104}
\check{c}_1(\kappa,\alpha)
=
\check{c}_{-1}(\kappa,\alpha)
=
\frac{\Gamma(\kappa+\mu+\tfrac12)
\Gamma(\kappa-\mu+\tfrac12)}
{2\pi}
\left(\frac{2e}{\kappa(4-\alpha^2)^{1/2}}\right)^{2\kappa}
\left(\frac{2-\alpha}{2+\alpha}\right)^{\mu}.
\end{equation}

We next determine the LG coefficients explicitly. In place of \cref{eq40}, introduce the variable $\check{\beta}$, which differs from $\beta$ only in the choice of branch. It has a branch cut along $[z_t,\check{z}_t]$ and is positive for $z\in(\check{z}_t,\infty)$, and thus
\begin{equation}
\label{eq105}
\check{\beta}
=
\left(\frac{z-z_t}{z-\check{z}_t}\right)^{1/2},
\end{equation}
so that $\check{\beta}\to1$ as $z\to\infty$ in any direction. Moreover, differentiating \cref{eq105} with respect to $\check{\xi}$ gives
\begin{equation}
\label{eq106}
\frac{d\check{\beta}}{d\check{\xi}}
=
-\frac{(1-\check{\beta}^{2})^{2}
\{\check{\beta}^{2}\check{z}_{t}-z_{t}\}}
{\check{\beta}^{2}(z_{t}-\check{z}_{t})^{2}}.
\end{equation}
It follows from \cref{eq34,eq35,eq36,eq38,eq106} that the LG coefficients in the present case are given by $E_s(\alpha,z)=\mathrm{E}_s(\alpha,-\check{\beta})$, where $\mathrm{E}_s(\alpha,\beta)$ are defined by \cref{eq47,eq48,eq50}. These are Laurent polynomials in $\check{\beta}$, and hence, from \cref{eq105}, each $(z-\check{z}_t)^{1/2}E_{2s+1}(\alpha,z)$ ($s=0,1,2,\ldots$) is meromorphic at $z=\check{z}_t$, as required in \cite{Dunster:2017:COA} for the validity of \cref{eq100,eq101} in $\check{Z}$.

Now from \cref{eq91,eq92} we deduce that
\begin{equation}
\label{eq107}
\lambda_{2s+1}^{(j)}(\alpha)
=-\lambda_{2s+1}(\alpha),
\quad
\lambda_{2s+2}^{(j)}(\alpha)=0
\quad (j=0,\pm 1,\, s=0,1,2,\ldots).
\end{equation}

Next, analogously to \cref{eq57,eq58}, let us define coefficients $\mathfrak{E}_{s}(\alpha,z)$ and $\tilde{\mathfrak{E}}_{s}(\alpha,z)$ by \cref{eq95,eq96} with $E_s(\alpha,z)=\mathrm{E}_s(\alpha,-\check{\beta})$; explicitly,
\begin{equation}
\label{eq108}
\mathfrak{E}_{s}(\alpha,z)
=
\mathrm{E}_s(\alpha,-\check{\beta})
+(-1)^s a_s s^{-1}\check{\xi}^{-s}
\end{equation}
and
\begin{equation}
\label{eq109}
\tilde{\mathfrak{E}}_{s}(\alpha,z)
=
\mathrm{E}_s(\alpha,-\check{\beta})
+(-1)^s \tilde a_s s^{-1}\check{\xi}^{-s}.
\end{equation}
Then, using \cref{eq10,eq21,eq97,eq98,eq99,eq102,eq103,eq104,eq107,eq108,eq109}, and recalling that $\lambda_{2s}(\alpha)=0$ ($s=1,2,3,\ldots$), we obtain
\begin{multline}
\label{eq110}
W_{\kappa,\mu}(\kappa z)
=\check{d}_{0}(\kappa,\alpha)
\left(\frac{z^2\check{\zeta}}
{\left(z-z_{t}\right)
\left(z-\check{z}_{t}\right)}\right)^{1/4}
\\ \times
\left\{
\mathrm{Ai}\left(\kappa^{2/3}\check{\zeta}\right)
\check{A}(\kappa,\alpha,z)
+
\mathrm{Ai}'\left(\kappa^{2/3}\check{\zeta}\right)
\check{B}(\kappa,\alpha,z)\right\}
\end{multline}
and
\begin{multline}
\label{eq111}
W_{-\kappa,\mu}(\kappa z e^{\mp \pi i})
=\check{d}_{\pm 1}(\kappa,\alpha)
\left(\frac{z^2\check{\zeta}}
{\left(z-z_{t}\right)
\left(z-\check{z}_{t}\right)}\right)^{1/4}
\\ \times
\left\{
\mathrm{Ai}_{\pm 1}\left(\kappa^{2/3}\check{\zeta}\right)
\check{A}(\kappa,\alpha,z)
+
\mathrm{Ai}'_{\pm 1}\left(\kappa^{2/3}\check{\zeta}\right)
\check{B}(\kappa,\alpha,z)\right\},
\end{multline}
where
\begin{equation}
\label{eq112}
\check{d}_{0}(\kappa,\alpha)
=
2\pi^{1/2}\kappa^{1/6}
\left(\frac{\kappa(4-\alpha^2)^{1/2}}{2e}\right)^\kappa
\left(\frac{2+\alpha}{2-\alpha}\right)^{\mu/2}
\end{equation}
and
\begin{equation}
\label{eq113}
\check{d}_{\pm 1}(\kappa,\alpha)
=
\frac{
4\pi^{3/2}\kappa^{1/6}
e^{\pm\kappa\pi i\mp\pi i/6}
}
{
\Gamma(\kappa+\mu+\tfrac12)
\Gamma(\kappa-\mu+\tfrac12)
}
\left(\frac{\kappa(4-\alpha^2)^{1/2}}{2e}\right)^\kappa
\left(\frac{2+\alpha}{2-\alpha}\right)^{\mu/2}.
\end{equation}
The coefficient functions then possess the asymptotic expansions
\begin{multline}
\label{eq114}
\check{A}(\kappa,\alpha,z)
\sim
\exp\left\{
-\sum\limits_{s=0}^{\infty}
\frac{\lambda_{2s+1}(\alpha)}{\kappa^{2s+1}}
\right\}
\exp\left\{\sum\limits_{s=1}^{\infty}
\frac{\tilde{\mathfrak{E}}_{2s}(\alpha,z)}{\kappa^{2s}}
\right\}
\\ \times
\cosh\left\{\sum\limits_{s=0}^{\infty}
\frac{\tilde{\mathfrak{E}}_{2s+1}(\alpha,z)}{\kappa^{2s+1}}
\right\}
\end{multline}
and
\begin{multline}
\label{eq115}
\check{B}(\kappa,\alpha,z)
\sim
\frac{1}
{\kappa^{1/3}\check{\zeta}^{1/2}}
\exp\left\{
-\sum\limits_{s=0}^{\infty}
\frac{\lambda_{2s+1}(\alpha)}{\kappa^{2s+1}}
\right\}
\exp\left\{\sum\limits_{s=1}^{\infty}
\frac{\mathfrak{E}_{2s}(\alpha,z)}{\kappa^{2s}}
\right\}
\\ \times
\sinh\left\{\sum\limits_{s=0}^{\infty}
\frac{\mathfrak{E}_{2s+1}(\alpha,z)}{\kappa^{2s+1}}
\right\},
\end{multline}
as $\kappa\to\infty$ for $z\in\check{Z}$ and $\alpha\in[0,\alpha_0]$.

We wish to extend these expansions to a full neighborhood of the turning point $z=\check{z}_t$, and to do so we proceed as in \cref{eq66,eq67,eq68,eq69}. Let $\tilde{\mathfrak{q}}_{1}(\alpha,z)=\tilde{\mathfrak{E}}_{1}(\alpha,z)$ and
\begin{equation}
\label{eq116}
\tilde{\mathfrak{q}}_{s}(\alpha,z)=\tilde{\mathfrak{E}}_{s}(\alpha,z)
+\frac{1}{s}\sum_{j=1}^{s-1}j
\tilde{\mathfrak{E}}_{j}(\alpha,z)\tilde{\mathfrak{q}}_{s-j}(\alpha,z)
\quad (s=2,3,4,\ldots)
\end{equation}
and let $\mathfrak{q}_{1}(\alpha,z)=\mathfrak{E}_{1}(\alpha,z)$, with
\begin{equation}
\label{eq117}
\mathfrak{q}_{s}(\alpha,z)=\mathfrak{E}_{s}(\alpha,z)
+\frac{1}{s}\sum_{j=1}^{s-1}j
\mathfrak{E}_{j}(\alpha,z)\mathfrak{q}_{s-j}(\alpha,z)
\quad (s=2,3,4,\ldots).
\end{equation}
Then, on setting $\check{\mathrm{A}}_{s}(\alpha,z)=\tilde{\mathfrak{q}}_{2s}(\alpha,z)$ ($s=1,2,3,\ldots$) and $\check{\mathrm{B}}_{s}(\alpha,z)=\check{\zeta}^{-1/2}\mathfrak{q}_{2s+1}(\alpha,z)$ ($s=0,1,2,\ldots$), we obtain
\begin{equation}
\label{eq118}
\check{A}(\kappa,\alpha,z) \sim 
\exp\left\{
-\sum\limits_{s=0}^{\infty}
\frac{\lambda_{2s+1}(\alpha)}{\kappa^{2s+1}}
\right\}
\left\{
1+\sum\limits_{s=1}^{\infty}
\frac{\check{\mathrm{A}}_{s}(\alpha,z)}{\kappa^{2s}}
\right\}
\end{equation}
and
\begin{equation}
\label{eq119}
\check{B}(\kappa,\alpha,z)
\sim
\frac{1}{\kappa^{4/3}}
\exp\left\{
-\sum\limits_{s=0}^{\infty}
\frac{\lambda_{2s+1}(\alpha)}{\kappa^{2s+1}}
\right\}
\sum\limits_{s=0}^{\infty}
\frac{\check{\mathrm{B}}_{s}(\alpha,z)}{\kappa^{2s}},
\end{equation}
as $\kappa\to\infty$, for $\alpha\in[0,\alpha_0]$ and $z\in\check{Z}\cup\{\check{z}_t\}$. The validity of these expansions at $z=\check{z}_t$ follows from \cite[Thm.~3.1]{Dunster:2021:NKF}. In particular, they hold uniformly for $\alpha\in[0,\alpha_0]$ and $z$ in the closed half-plane $\Re(z)\ge z_t+\delta_1$, where $\delta_1\in(0,\check{z}_t-z_t)$ is an arbitrary fixed constant.

\subsection{Main Airy expansions and numerical results}

We now state the main results of this section. Let $\alpha$, $z_t$, $\check{z}_t$, $\check{\zeta}$, and $\check{\xi}$ be given by \cref{eq19,eq22,eq82}, and let $\check{\beta}$ be given by \cref{eq105}. Let $\mathrm{E}_s(\alpha,\beta)$ and $\lambda_{2s+1}(\alpha)$ be defined by \cref{eq45,eq47,eq48,eq50}. Define the sequences $\{a_s\}$ and $\{\tilde a_s\}$ by \cref{eq94}, and $\mathfrak{E}_s(\alpha,z)$ and $\tilde{\mathfrak{E}}_s(\alpha,z)$ by \cref{eq108,eq109}. Finally, define $\tilde{\mathfrak{q}}_s(\alpha,z)$ and $\mathfrak{q}_s(\alpha,z)$ by \cref{eq116,eq117}, and hence $\check{\mathrm{A}}_s(\alpha,z)$ and $\check{\mathrm{B}}_s(\alpha,z)$ as stated immediately following \cref{eq117}. The results below are then obtained from \cref{eq110,eq112} for \cref{eq120}, and from \cref{eq11}, \cref{eq12}, \cref{eq111,eq112,eq113} and \cite[Eq.~9.2.11]{NIST:DLMF} for \cref{eq121}.
\begin{theorem}
\label{thm:Airy}
\begin{multline}
\label{eq120}
W_{\kappa,\mu}(\kappa z)
=2\pi^{1/2}\kappa^{1/6}
\left(\frac{\kappa(4-\alpha^2)^{1/2}}{2e}\right)^\kappa
\left(\frac{2+\alpha}{2-\alpha}\right)^{\mu/2}
\left(\frac{z^2\check{\zeta}}
{\left(z-z_{t}\right)
\left(z-\check{z}_{t}\right)}\right)^{1/4}
\\ \times
\left\{
\mathrm{Ai}\left(\kappa^{2/3}\check{\zeta}\right)
\check{A}(\kappa,\alpha,z)
+
\mathrm{Ai}'\left(\kappa^{2/3}\check{\zeta}\right)
\check{B}(\kappa,\alpha,z)\right\}
\end{multline}
and 
\begin{multline}
\label{eq121}
\mathbf{M}_{\kappa,\pm\mu}(\kappa z)
=\frac{2\pi^{1/2}\kappa^{1/6}}
{\Gamma(\kappa\pm\mu+\tfrac12)}
\left(\frac{\kappa(4-\alpha^2)^{1/2}}{2e}\right)^\kappa
\left(\frac{2+\alpha}{2-\alpha}\right)^{\mu/2}
\left(\frac{z^2\check{\zeta}}
{\left(z-z_{t}\right)
\left(z-\check{z}_{t}\right)}\right)^{1/4}
\\ \times
\left\{
\mathscr{A}\left(\kappa^{2/3}\check{\zeta},
(\kappa\mp\mu-\tfrac12)\pi\right)
\check{A}(\kappa,\alpha,z)
+
\mathscr{A}'\left(\kappa^{2/3}\check{\zeta},
(\kappa\mp\mu-\tfrac12)\pi\right)
\check{B}(\kappa,\alpha,z)
\right\},
\end{multline}
where
\begin{equation}
\label{eq122}
\mathscr{A}(z,\theta)
=\cos(\theta)\mathrm{Ai}(z)-\sin(\theta)\mathrm{Bi}(z)
\end{equation}
and $\mathscr{A}'(z,\theta)=\partial\mathscr{A}(z,\theta)/\partial z$. The functions $\check{A}(\kappa,\alpha,z)$ and $\check{B}(\kappa,\alpha,z)$ are analytic at $z=\check{z}_t$, and possess the asymptotic expansions \cref{eq118,eq119} as $\kappa\to\infty$, uniformly for $0\le\mu\le(1-\delta_0)\kappa<\kappa$, $\Re(z)\ge z_t+\delta_1$, and $-\frac12 \pi<\arg(z)<\frac12\pi$, where $\delta_0\in(0,1)$ and $\delta_1\in(0,\check{z}_t-z_t)$ are arbitrary fixed constants.

\end{theorem}

\begin{remark}
As in \cref{thm:Bessel}, the region of validity in the $z$-plane is larger than that stated in \cref{thm:Airy}. For simplicity, we again restrict attention to a convenient half-plane subregion which contains $z=\check{z}_t$. Moreover, the constants $\delta_1$ may be chosen so that the two half-plane regions of \cref{thm:Bessel,thm:Airy} overlap and hence cover the entire principal $z$-plane.
\end{remark}

We next illustrate the accuracy of the Airy expansion for $W_{\kappa,\mu}(\kappa z)$ for real $z>z_t$. As in the Bessel case, a suitable envelope is needed to measure the relative error through the oscillatory region, and to do so we first introduce a numerically satisfactory companion solution. Motivated by \cref{eq10}, define
\begin{multline}
\label{eq123}
\check{W}_{\kappa,\mu}(t)
=\frac{\Gamma(\kappa+\mu+\tfrac12)
\Gamma(\kappa-\mu+\tfrac12)}{2\pi}
\\ \times 
\left\{
e^{-\kappa\pi i}W_{-\kappa,\mu}(t e^{-\pi i})
+e^{\kappa\pi i}W_{-\kappa,\mu}(t e^{\pi i})
\right\}.
\end{multline}
It then follows from \cref{eq111,eq123} and the Airy connection formulas \cite[Eq.~9.2.11]{NIST:DLMF} that
\begin{multline}
\label{eq124}
\check{W}_{\kappa,\mu}(\kappa z)
=
\check d_0(\kappa,\alpha)
\left(\frac{z^2\check\zeta}
{(z-z_t)(z-\check z_t)}\right)^{1/4}
\\ \times
\left\{
\mathrm{Bi}\left(\kappa^{2/3}\check\zeta\right)\check A(\kappa,\alpha,z)
+
\mathrm{Bi}'\left(\kappa^{2/3}\check\zeta\right)\check B(\kappa,\alpha,z)
\right\}.
\end{multline}
Thus $\check{W}_{\kappa,\mu}(\kappa z)$ provides the natural companion involving $\mathrm{Bi}$ to the $\mathrm{Ai}$ representation of $W_{\kappa,\mu}(\kappa z)$ in \cref{eq120}. We accordingly define the positive envelope
\begin{equation}
\label{eq125}
\check{\mathsf{M}}(\kappa,\mu,z)=
    \begin{cases}
        \left[
        \left\{W_{\kappa,\mu}(\kappa z)\right\}^{2}
        +
        \left\{\check{W}_{\kappa,\mu}(\kappa z)\right\}^{2}
        \right]^{1/2}
        & (z_t < z \leq \check \omega_{\kappa,\mu}),\\[5pt]
        W_{\kappa,\mu}(\kappa z)
        & (\check \omega_{\kappa,\mu} < z<\infty),
    \end{cases}
\end{equation}
where $\check \omega_{\kappa,\mu}$ is the largest zero, in the $z$-variable, of $\check{W}_{\kappa,\mu}(\kappa z)$. On $(z_t,\check\omega_{\kappa,\mu}]$ the first expression closely follows the amplitude of $W_{\kappa,\mu}(\kappa z)$ in the oscillatory region, whereas beyond $\check\omega_{\kappa,\mu}$ the envelope coincides with $W_{\kappa,\mu}(\kappa z)$.

\begin{figure}[!htb]
 \centering
 \includegraphics[
 width=0.9\textwidth,keepaspectratio]{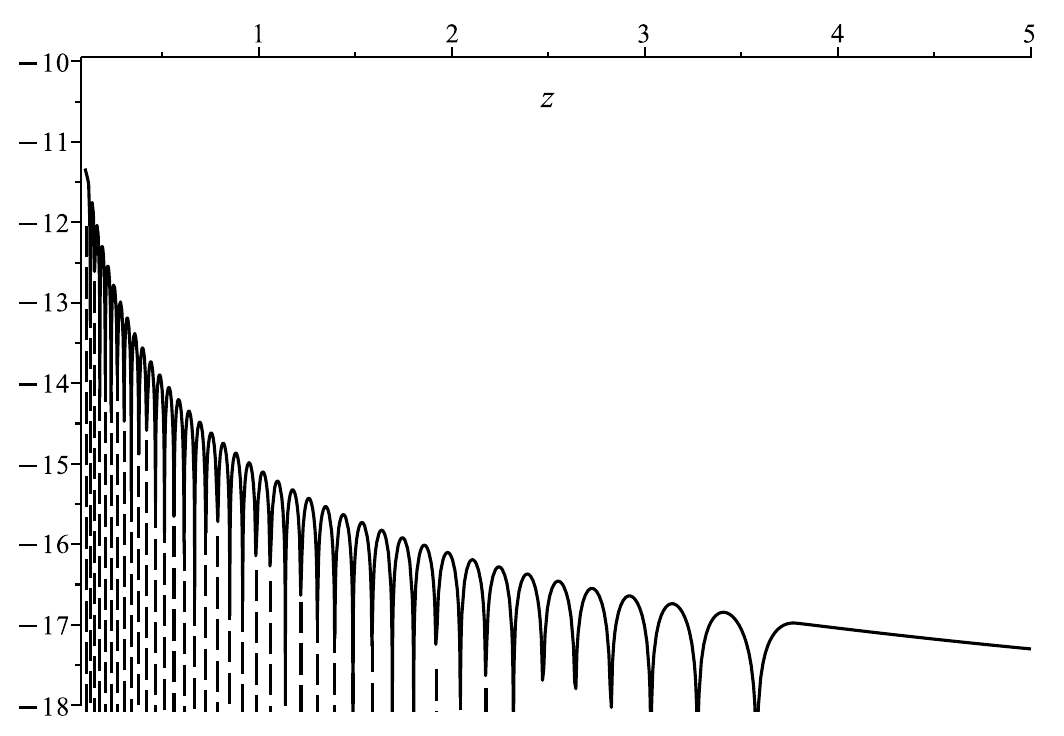}
\caption{Graph of $\check{\Omega}_{4}(\kappa,\mu,z)$ for $\kappa=50$, $\alpha=0.01$ ($\mu=0.25$)}
 \label{fig:hatomega0.01}
\end{figure}

\begin{figure}[!htb]
 \centering
 \includegraphics[
 width=0.9\textwidth,keepaspectratio]{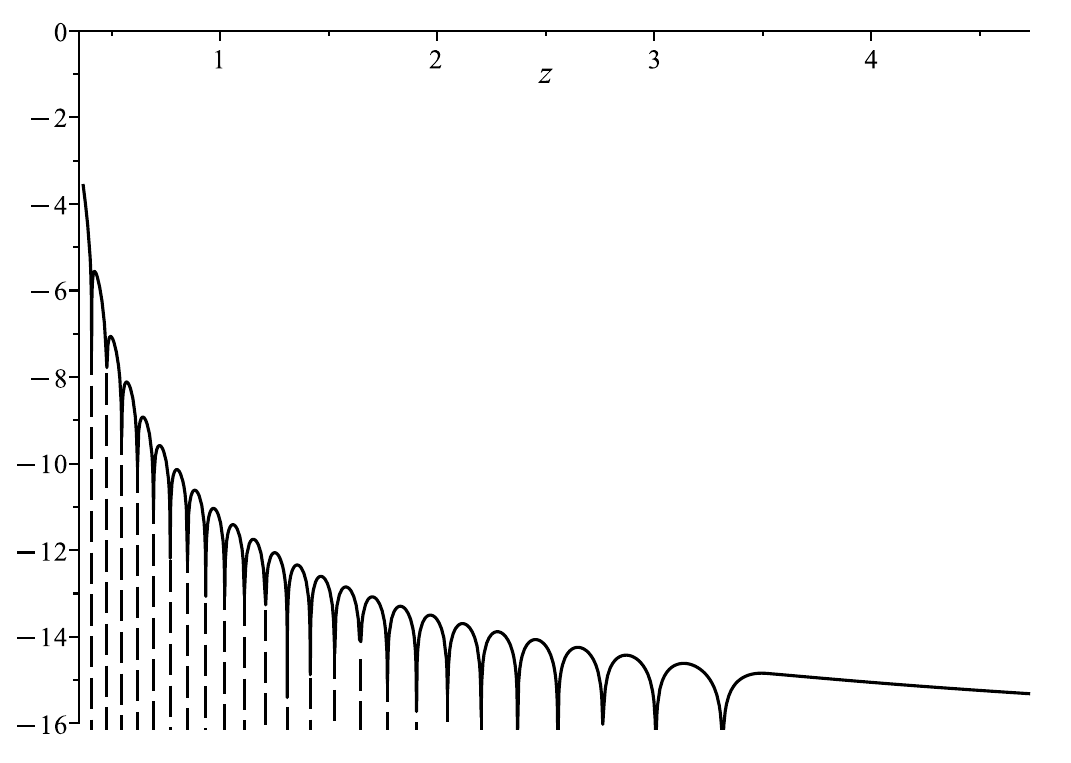}
 \caption{Graph of $\check{\Omega}_{4}(\kappa,\mu,z)$ for $\kappa=50$, $\alpha=1$ ($\mu=25$)}
 \label{fig:hatomega1}
\end{figure}

\begin{figure}[!htb]
 \centering
 \includegraphics[
 width=0.9\textwidth,keepaspectratio]{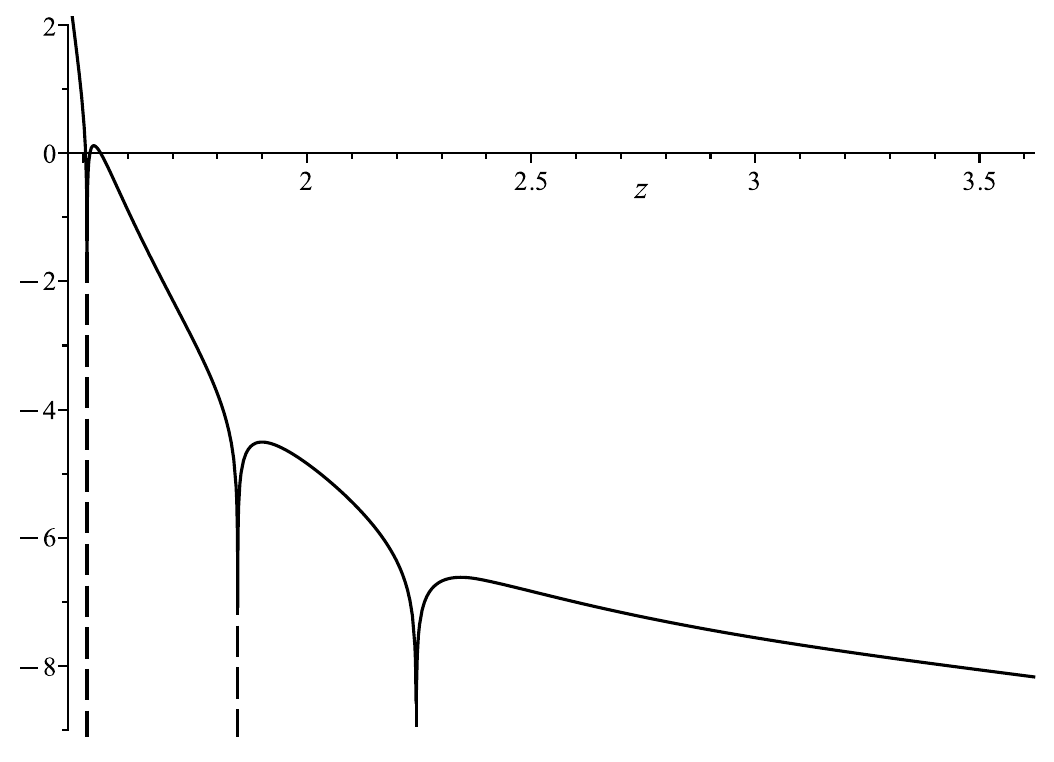}
\caption{Graph of $\check{\Omega}_{4}(\kappa,\mu,z)$ for $\kappa=50$, $\alpha=1.9$ ($\mu=47.5$)}
 \label{fig:hatomega1.9}
\end{figure}

For comparison with the exact solution $W_{\kappa,\mu}(\kappa z)$, define $W_4(\kappa,\mu,z)$ for $z>z_t$ by retaining the first four terms of the expansions \cref{eq118,eq119}; thus,
\begin{multline}
\label{eq126}
W_4(\kappa,\mu,z)
=
\check d_0(\kappa,\alpha)
\left(\frac{z^2\check\zeta}
{(z-z_t)(z-\check z_t)}\right)^{1/4}
\exp\left\{
-\sum\limits_{s=0}^{3}
\frac{\lambda_{2s+1}(\alpha)}{\kappa^{2s+1}}
\right\}
\\ \times
\left[
\mathrm{Ai}\left(\kappa^{2/3}\check\zeta\right)
\left\{
1+\sum\limits_{s=1}^{3}
\frac{\check{\mathrm A}_{s}(\alpha,z)}{\kappa^{2s}}
\right\}
+
\frac{\mathrm{Ai}'\left(\kappa^{2/3}\check\zeta\right)}
{\kappa^{4/3}}
\sum\limits_{s=0}^{3}
\frac{\check{\mathrm B}_{s}(\alpha,z)}{\kappa^{2s}}
\right].
\end{multline}
For numerical computation near the turning point $z=\check{z}_t$, we handled the apparent singularities in $\check{\mathrm A}_{s}(\alpha,z)$ and $\check{\mathrm B}_{s}(\alpha,z)$ by using their Taylor expansions about $z=\check{z}_t$.

Similarly to \cref{eq81}, the relative error can now be measured by
\begin{equation}
\label{eq127}
\check{\Omega}_{4}(\kappa,\mu,z)
=
\log_{10}\left\{
\frac{
\left|W_{\kappa,\mu}(\kappa z)-W_4(\kappa,\mu,z)\right|
}
{\check{\mathsf{M}}(\kappa,\mu,z)}
\right\}.
\end{equation}
For $\kappa=50$, \cref{fig:hatomega0.01,fig:hatomega1,fig:hatomega1.9} show $\check{\Omega}_4(\kappa,\mu,z)$ for $\alpha=0.01$, $1$, and $1.9$, respectively. The behavior is complementary to that observed for the Bessel expansions. Here the approximation becomes increasingly accurate as $z$ moves away from the first turning point $z_t$, with no loss of precision on passing through the second turning point $z=\check z_t$. The accuracy continues to improve into the unbounded interval $z>\check z_t$, in accordance with the uniformity of the Airy expansion there. For $\alpha=0.01$ and $1$, the numerical results show very high precision in a broad region containing $\check z_t$ and extending to its right. For $\alpha=1.9$, where the two turning points are appreciably closer, the approximation is initially less accurate near $z_t$, but improves rapidly as $z$ increases and remains good through and beyond $\check z_t$. Thus the numerical results for the Bessel and Airy expansions exhibit the complementary behavior expected from their respective regions of validity.

\section{Generalized Laguerre polynomials}
\label{sec:Laguerre}

As a consequence of the preceding results, we obtain uniform asymptotic expansions for the generalized Laguerre polynomials. These apply when
\begin{equation}
\label{eq128}
\kappa=n+\mu+\tfrac12
\quad
(n=0,1,2,\ldots),
\end{equation}
and for this case, on referring to \cref{eq13} and the standard relation between the Whittaker function and the generalized Laguerre polynomial \cite[Eq.~13.18.17]{NIST:DLMF}, we have the identity
\begin{equation}
\label{eq129}
L_n^{(2\mu)}(\kappa z)
=
\frac{(-1)^n e^{\kappa z/2}}
{n!(\kappa z)^{\mu+1/2}}
W_{\kappa,|\mu|}(\kappa z).
\end{equation}

Assume first that $\mu \ge 0$. Then since $(\mu-\kappa+\tfrac12)\pi=-n\pi$, it follows from \cref{eq76} that $\mathscr{C}_{2\mu}(x,(\mu-\kappa+\tfrac12)\pi)=(-1)^nJ_{2\mu}(x)$. Hence, on substituting \cref{eq75} into \cref{eq129} and using $\Gamma(\kappa-\mu+\tfrac12)=n!$, we obtain the Bessel expansion
\begin{multline}
\label{eq130}
L_n^{(2\mu)}(\kappa z)
=
\frac{\Gamma(\kappa+\mu+\tfrac12)}
{\pi\kappa^{\mu+1/2}}
e^{\kappa z/2}z^{-\mu}
\left(\frac{2e}{\kappa(4-\alpha^2)^{1/2}}\right)^\kappa
\left(\frac{2-\alpha}{2+\alpha}\right)^{\mu/2}
\\
\times
\left(\frac{\alpha^2-\zeta}
{(z_t-z)(\check z_t-z)}\right)^{1/4}
\Biggl\{
J_{2\mu}\left(\kappa\zeta^{1/2}\right)
A_\nu(\kappa,\alpha,z)
+
\frac{\zeta}{\kappa}
\frac{\partial J_{2\mu}\left(\kappa\zeta^{1/2}\right)}
{\partial\zeta}
B_\nu(\kappa,\alpha,z)
\Biggr\}.
\end{multline}

Likewise, substitution of the Airy expansion \cref{eq110} into \cref{eq129} gives
\begin{multline}
\label{eq131}
L_n^{(2\mu)}(\kappa z)
=(-1)^n
\frac{\check d_0(\kappa,\alpha)}
{n!\kappa^{\mu+1/2}}
e^{\kappa z/2}z^{-\mu}
\left(\frac{\check\zeta}
{(z-z_t)(z-\check z_t)}\right)^{1/4}
\\
\times
\left\{
\mathrm{Ai}\left(\kappa^{2/3}\check\zeta\right)
\check A(\kappa,\alpha,z)
+
\mathrm{Ai}'\left(\kappa^{2/3}\check\zeta\right)
\check B(\kappa,\alpha,z)
\right\}.
\end{multline}
Here $\alpha$, $\zeta$, $\check{\zeta}$, $z_t$, $\check z_t$, and all the coefficient functions are as defined in \cref{sec:Bessel,sec:Airy}. The Bessel and Airy expansions \cref{eq68,eq69,eq118,eq119} hold for \cref{eq130,eq131} and inherit the respective regions of validity and uniformity properties established there. Thus, these asymptotic expansions for \cref{eq130} hold uniformly for $\Re(z)\le\check{z}_t-\delta_1$, while those used in \cref{eq131} hold uniformly for $\Re(z)\ge z_t+\delta_1$, where $\delta_1\in(0,\check{z}_t-z_t)$ is an arbitrary fixed constant. Moreover, from \cref{eq25,eq128}, in both cases these expansions hold uniformly for
\begin{equation}
\label{eq132}
0\le\mu\le\frac{1-\delta_0}{\delta_0}\left(n+\tfrac12\right),
\end{equation}
as $n\to\infty$, where $\delta_0\in(0,1)$. Thus the Laguerre parameter $2\mu$ may be larger than the degree $n$; as $\delta_0\to0$, the admissible ratio $\mu/n$ is unbounded.

Next consider $\mu<0$. It turns out that \cref{eq130,eq131}, together with their associated asymptotic expansions \cref{eq68,eq69,eq118,eq119}, remain valid in this case, provided
\begin{equation}
\label{eq133}
-\frac{1-\delta_0}{2-\delta_0}\left(n+\tfrac12\right)
\le \mu<0,
\end{equation}
again where $\delta_0\in(0,1)$. To see this, we can again use \cref{eq129}, with \cref{eq75} applied after replacing $\mu$ by $|\mu|$ and $\alpha$ by $|\alpha|$. In order to do so, we require $\kappa>0$ and the condition \cref{eq25} with $\mu$ replaced by $|\mu|$. Both requirements are satisfied under \cref{eq133}; indeed, \cref{eq133} is equivalent to $|\mu|\le(1-\delta_0)\kappa$, and also implies $\kappa=n+\mu+\tfrac12>0$. Now, since $\kappa=n+\mu+\tfrac12=n-|\mu|+\tfrac12$, we have $(|\mu|-\kappa+\tfrac12)\pi=(2|\mu|-n)\pi$. Hence, from \cref{eq76} and \cite[Eq.~10.4.5]{NIST:DLMF},
\begin{multline}
\label{eq134}
\mathscr{C}_{2|\mu|}\left(\kappa \zeta^{1/2},(|\mu|-\kappa+\tfrac12)\pi\right)
=\mathscr{C}_{2|\mu|}\left(\kappa \zeta^{1/2},(2|\mu|-n)\pi\right)
\\
=(-1)^n\left\{\cos\left(2|\mu|\pi\right)J_{2|\mu|}\left(\kappa \zeta^{1/2}\right)
-\sin\left(2|\mu|\pi\right)Y_{2|\mu|}\left(\kappa \zeta^{1/2}\right)\right\}
\\
=(-1)^nJ_{-2|\mu|}\left(\kappa \zeta^{1/2}\right)
=(-1)^nJ_{2\mu}\left(\kappa \zeta^{1/2}\right),
\end{multline}
which is the same as the positive $\mu$ case. The turning points $z_t$ and $\check{z}_t$ are also unchanged, since they depend on $\alpha$ only through $\alpha^2$, and consequently the associated variables $\xi$ and $\zeta$ are unchanged. The coefficient functions $A_\nu(\kappa,\alpha,z)$ and $B_\nu(\kappa,\alpha,z)$ are likewise unchanged under $\alpha\mapsto|\alpha|$. Thus, using \cref{eq134} in \cref{eq75} and substituting the result into \cref{eq129}, we recover exactly the same Bessel-type expansion \cref{eq130} for $\mu<0$.

The same argument applies to the Airy-type expansions: replacing $\mu$ by $|\mu|$ and $\alpha$ by $|\alpha|$ leaves the turning points, the Airy variable, and the corresponding coefficient functions unchanged, so the resulting expansions retain exactly the same form.

In summary, \cref{eq130,eq131}, when used in combination, provide uniform asymptotic expansions as $n\to\infty$ throughout the unbounded complex $z$-plane, including at $z=0$, uniformly for $\mu$ in the ranges \cref{eq132,eq133}. Most notably, the Laguerre parameter $a=2\mu$ ranges over the large interval
$-(1-\delta)n\le a\le\Delta n$, where $\delta$ and $\Delta$ are arbitrary fixed positive constants.

We finally remark that in \cite{Dunster:2018:USE}, similar Airy-type expansions were obtained for the Laguerre polynomials in neighborhoods of the turning points $z=z_t$ and $z=\check{z}_t$, which, as here, remain bounded away from each other, with separate expansions constructed at each point. The expansion about $z=\check{z}_t$ is analogous to \cref{eq131}. However, the expansion of \cite{Dunster:2018:USE} valid at $z=z_t$ requires this turning point to remain bounded away from the pole at $z=0$, and hence $\alpha$ must be bounded away from zero. Thus the parameter range, although large, is more restricted than in the present results.

\section*{Acknowledgement}
Financial support from Ministerio de Ciencia e Innovación project PID2024-159583NB-I00 (MICIU/ AEI / 10.13039/501100011033 / FEDER, UE) is acknowledged.

\makeatletter
\interlinepenalty=10000

\bibliographystyle{siamplain}
\bibliography{biblio}

@article{Dunster:2025:LCT,
  title={Asymptotic expansions for solutions of differential equations having coalescing turning points, with an application to {L}egendre functions},
  author={Dunster, T. M.},
  journal={Stud. Appl. Math.},
  year={2025},
  note={},
  eprint={},
  doi={10.1111/sapm.70138},
  archivePrefix={},
  primaryClass={math.CA},
  url={}
}

@article{Dunster:2026:TPD,
  title={Asymptotic expansions for solutions of differential equations having a coalescing turning point and double pole, with an application to {L}egendre functions},
  author={Dunster, T. M.},
  journal={Stud. Appl. Math.},
  volume={156},
  number={5},
  pages={e70240},
  year={2026},
  note={},
  eprint={},
  doi={10.1111/sapm.70240},
  archivePrefix={},
  primaryClass={math.CA},
  url={}
}

@article{Temme:1978:UAC,
    author = {Temme, N. M.},
    title = {Uniform asymptotic expansions of confluent hypergeometric functions},
    journal = {IMA J. Appl. Math.},
    volume = {22},
    year = {1978},
    number = {2},
    pages = {215-223},
    doi = {10.1093/imamat/22.2.215}
}

@article{Dunster:2018:USE,
	Author = {Dunster, T. M. and Gil, A. and Segura, J.},
	Day = {01},
	Doi = {10.1007/s10444-018-9589-5},
	Issn = {1572-9044},
	Journal = {Adv. Comput. Math.},
	Month = {Oct},
	Number = {5},
	Pages = {1441-1474},
	Title = {Uniform asymptotic expansions for {L}aguerre polynomials and related confluent hypergeometric functions},
	Url = {},
	Volume = {44},
	Year = {2018}
}

@article{Dunster:2021:SEB,
	Author = {Dunster, T. M. and Gil, A. and Segura, J.},
	Journal = {Anal. Appl.},
	Doi ={10.1142/S0219530520500104},
	Title = {Simplified error bounds for turning point expansions},
	Year = {2021},
	Number = {4},
	Volume = {19},
	Pages = {647-678}}

@article{Dunster:2020:LGE,
	Author = {Dunster, T. M.},
	Doi = {10.1017/prm.2018.117},
	Journal = {Proc. Roy. Soc. Edinburgh Sec. A},
	Number = {3},
	Pages = {1289-1311},
	Publisher = {Royal Society of Edinburgh Scotland Foundation},
	Title = {Liouville-{G}reen expansions of exponential form, with an application to modified {B}essel functions},
	Volume = {150},
	Year = {2020}}

@misc{NIST:DLMF,
	howpublished = {Release 1.1.6 of 2022-06-30},
	key = {{\relax DLMF}},
	note = {F.~W.~J. Olver, A.~B. {Olde Daalhuis}, D.~W. Lozier, B.~I. Schneider, R.~F. Boisvert, C.~W. Clark, B.~R. Miller, B.~V. Saunders, H.~S. Cohl, and M.~A. McClain, eds.},
	title = {{\it NIST Digital Library of Mathematical Functions}},
	url = {http://dlmf.nist.gov/}}

@article{Dunster:2021:UAW,
    author = {Dunster, T. M.},
    title = {Uniform asymptotic expansions for the {W}hittaker functions ${M}_{\kappa,\mu}(z)$ and ${W}_{\kappa,\mu}(z)$ with $\mu$ large},
    journal = {Proc. R. Soc. A},
    volume = {477},
    year = {2021},
    number = {2252},
    pages = {20210360},
    doi = {10.1098/rspa.2021.0360}
}

@article{Dunster:2025:SAR,
  author       = {T. M. Dunster},
  title        = {Simplified {A}iry function Asymptotic expansions for Reverse Generalised {B}essel Polynomials},
  year         = {2026},
  number         = {1},
  volume         = {28},
  journal      = {J. Classical Anal.},
  Doi={10.7153/jca-2026-28-01},
  URL = {},
  note         = {}}

@book{Buchholz:1969:CHF,
	Author = {Buchholz, H.},
	Publisher = {Springer, Berlin, Heidelberg},
	Title = {The confluent hypergeometric function},
	Doi = {10.1007/978-3-642-88396-5},
	Year = {1969}}

@book{Olver:1997:ASF,
	Address = {Wellesley, MA},
	Author = {Olver, F. W. J.},
	Isbn = {1-56881-069-5},
	Mrclass = {41-02 (33Cxx 41A60 65D20)},
	Mrnumber = {MR1429619 (97i:41001)},
	Note = {Reprint of the 1974 original [Academic Press, New York]},
	Pages = {xviii+572},
	Publisher = {A K Peters Ltd.},
	Series = {AKP Classics},
	Title = {Asymptotics and special functions},
	Year = {1997}}

@article{Dunster:2017:COA,
	Author = {Dunster, T. M. and Gil, A. and Segura, J.},
	Doi = {10.1007/s00365-017-9372-8},
	Journal = {Constr. Approx.},
	Number = {3},
	Pages = {645-675},
	Title = {Computation of asymptotic expansions of turning point problems via {C}auchy's integral formula: Bessel functions},
	Volume = {46},
	Year = {2017}}

@article{Dunster:1989:UAE,
	Author = {Dunster, T. M.},
	Doi = {10.1137/0520052},
	Fjournal = {SIAM Journal on Mathematical Analysis},
	Issn = {0036-1410},
	Journal = {SIAM J. Math. Anal.},
	Mrclass = {33A30 (34E20)},
	Mrnumber = {990876},
	Mrreviewer = {F. W. J. Olver},
	Number = {3},
	Pages = {744-760},
	Title = {Uniform asymptotic expansions for {W}hittaker's confluent hypergeometric functions},
	Url = {},
	Volume = {20},
	Year = {1989}}

@article{Olver:1980:WFW,
	Author = {Olver, F. W. J.},
	Doi = {10.1017/S0308210500012130},
	Fjournal = {Proceedings of the Royal Society of Edinburgh. Section A. Mathematical and Physical Sciences},
	Issn = {0308-2105},
	Journal = {Proc. Roy. Soc. Edinburgh Sect. A},
	Mrclass = {34E20 (33A40)},
	Mrnumber = {592550},
	Mrreviewer = {P. F. Hsieh},
	Number = {3-4},
	Pages = {213-234},
	Title = {Whittaker functions with both parameters large: uniform approximations in terms of parabolic cylinder functions},
	Url = {http://dx.doi.org/10.1017/S0308210500012130},
	Volume = {86},
	Year = {1980}}

@article{Dunster:2021:NKF,
	Author = {Dunster, T. M.},
	Doi = {10.1137/21M1401590},
	Fjournal = {SIAM Journal on Mathematical Analysis},
	Journal = {SIAM J. Math. Anal.},
	Number = {5},
	Pages = {5915-5947},
	Title = {Nield-{K}uznetsov functions and {L}aplace transforms of parabolic cylinder functions},
	Url ={},
	Volume = {53},
	Year = {2021}}

@book{Slater:1960:CHF,
    author = {Slater, L. J.},
    title = {Confluent hypergeometric functions},
    publisher = {Cambridge University Press},
    address = {New York},
    year = {1960},
    pages = {xi+247},
    doi = {}
}

@article{Kubota:1969:SES,
    author = {Kubota, Nobuko and Yamaguchi, Mitsuko and Iwata, Giiti},
    title = {{S}chroedinger equations soluble in terms of confluent hypergeometric functions},
    journal = {Natur. Sci. Rep. Ochanomizu Univ.},
    volume = {20},
    year = {1969},
    number = {2},
    pages = {39-44},
    doi = {},
    url = {https://teapot.lib.ocha.ac.jp/records/34941}
}

@article{Stovicek:2024:CGW,
    author = {{\v S}{\v t}ov{\'i}{\v c}ek, Pavel},
    title = {{C}oulomb {G}reen's function and an addition formula for the {W}hittaker functions},
    journal = {J. Math. Phys.},
    volume = {65},
    year = {2024},
    number = {2},
    pages = {023503},
    doi = {10.1063/5.0184924}
}

@article{Heading:1965:RIP,
    author = {Heading, John},
    title = {Refractive index profiles based on the hypergeometric equation and the confluent hypergeometric equation},
    journal = {Proc. Cambridge Philos. Soc.},
    volume = {61},
    year = {1965},
    number = {4},
    pages = {897-913},
    doi = {10.1017/S0305004100039293}
}

@article{Dyson:1960:SIA,
    author = {Dyson, Freeman J.},
    title = {Stability of an idealized atmosphere. {II}. Zeros of the confluent hypergeometric function},
    journal = {Phys. Fluids},
    volume = {3},
    year = {1960},
    number = {2},
    pages = {155-157},
    doi = {10.1063/1.1706012}
}

@article{Saad:2003:ICH,
    author = {Saad, Nasser and Hall, Richard L.},
    title = {Integrals containing confluent hypergeometric functions with applications to perturbed singular potentials},
    journal = {J. Phys. A},
    volume = {36},
    year = {2003},
    number = {28},
    pages = {7771-7788},
    doi = {10.1088/0305-4470/36/28/307}
}

@article{Nogueira:2016:LCB,
    author = {Nogueira, P. H. F. and de Castro, A. S. and Pimentel, D. R. M.},
    title = {A large class of bound-state solutions of the {S}chr\"odinger equation via {L}aplace transform of the confluent hypergeometric equation},
    journal = {J. Math. Chem.},
    volume = {54},
    year = {2016},
    number = {6},
    pages = {1287-1295},
    doi = {10.1007/s10910-016-0621-z}
}

@article{Esparza:1999:ZWF,
    author = {Esparza, Javier and L{\'o}pez, Jos{\'e} L. and Sesma, Javier},
    title = {Zeros of the {W}hittaker function associated to {C}oulomb waves},
    journal = {IMA J. Appl. Math.},
    volume = {63},
    year = {1999},
    number = {1},
    pages = {71-87},
    doi = {10.1093/imamat/63.1.71}
}

@article{Kazarinoff:1955:AEW,
    author = {Kazarinoff, Nicholas D.},
    title = {Asymptotic expansions for the {W}hittaker functions of large complex order $m$},
    journal = {Trans. Amer. Math. Soc.},
    volume = {78},
    year = {1955},
    number = {2},
    pages = {305-328},
    doi = {10.2307/1993065}
}

@article{Kazarinoff:1957:AFW,
    author = {Kazarinoff, Nicholas D.},
    title = {Asymptotic forms for the {W}hittaker functions with both parameters large},
    journal = {J. Math. Mech.},
    volume = {6},
    year = {1957},
    number = {3},
    pages = {341-360},
    doi = {},
    url = {http://www.jstor.org/stable/24900469}
}

@book{Skovgaard:1966:UAE,
    author = {Skovgaard, Helge},
    title = {Uniform asymptotic expansions of confluent hypergeometric functions and {W}hittaker functions},
    publisher = {Jul. Gjellerups Forlag},
    address = {Copenhagen},
    year = {1966},
    pages = {iv+92},
    doi = {}
}
\end{document}